\documentclass[11pt]{amsart}
\usepackage{latexsym}
\usepackage{amsfonts}

\usepackage{amsmath,amsthm,amssymb}

\newtheorem{thm}{Theorem}[section]
\newtheorem{lem}[thm]{Lemma}
\newtheorem{cor}[thm]{Corollary}

\newtheorem{prop}[thm]{Proposition}
\theoremstyle{definition}
\newtheorem{defn}[thm]{Definition}
\newtheorem{notation}[thm]{Notation}
\newtheorem{conv}[thm]{Convention}
\newtheorem{rem}[thm]{Remark}

\newtheorem{propdfn}[thm]{Proposition-Definition}
\newtheorem{lede}[thm]{Lemma-Definition}

\begin{document}

\author[Ilya Kapovich]{Ilya Kapovich}

\address{\tt Department of Mathematics, University of Illinois at
  Urbana-Champaign, 1409 West Green Street, Urbana, IL 61801, USA
  \newline http://www.math.uiuc.edu/\~{}kapovich/} \email{\tt
  kapovich@math.uiuc.edu}

\thanks{The author was supported by the NSF grant DMS\#0404991 and by
  Centre de Recerca Matem\`atica at Barcelona}

\title[Currents on free groups]{Currents on free groups}

\begin{abstract}
 We study the properties of geodesic currents on free groups,
 particularly the ``intersection form" that is similar to Bonahon's
 notion of the intersection number between geodesic currents on
 hyperbolic surfaces.
\end{abstract}



\maketitle


\section{Introduction}\label{intro}

A \emph{geodesic current} on a word-hyperbolic group $G$ is a positive
$G$-invariant Borel measure on the space $\partial^2 G:=\{(x,y):
x,y\in \partial G, x\ne y\}$.  The study of geodesic currents on free
groups is motivated by obtaining information about the geometry and
dynamics of individual automorphisms as well as of groups of
automorphisms of a free group.  A similar programme has proved to be
very successful for the case of surface groups and hyperbolic
surfaces.  There Bonahon's foundational work~\cite{Bo86,Bo88} showed
the relevance of currents to the study of the geometry of the
Teichmuller space and the Thurston compactification of it, and to
understanding the dynamical properties of the mapping class group and
its individual elements.  Other interesting and important results
about geodesic currents in the hyperbolic surface case can be found in
\cite{Bridg,BT,Ha,Sa,Pi} and other sources.

We believe that in the free group case the study of currents is
particularly promising, in part since they are naturally defined in
the context of symbolic dynamics which is, in a sense, ``native'' to
the free group case. Some examples of geometric information about free
group automorphisms obtained by studying currents are contained in the
previous work of the author~\cite{Ka} as well as in \cite{KKS}.  Two
very important precursors and sources of inspiration for the present
paper are the article of Bonahon~\cite{Bo91}, where currents on
general word-hyperbolic groups are considered, and the 1995 doctoral
dissertation of Reiner Martin from UCLA~\cite{Ma}.

The main purpose of this paper is to collect together and clarify
various background facts and ideas related to geodesic currents on
free groups as well as to explain the relationships between them. We
also aim to set up the basic machinery (and even notations) for future
use and to clarify a number of typically confusing points (such as
those related to left and right actions of $Out(F_n)$ on the outer
space and on the space of currents.  We should make it clear that many
basic results here are not new and, in many instances are already
present, either implicitly or explicitly, in Bonahon and Martin's work
as well as in the author's article~\cite{Ka}.

Although the present paper is broad in scope, its main underlying
theme is to explain the significance and the geometric meaning of the
``intersection form'', which we believe to be a fundamentally
important object.

A central point of Bonahon's work~\cite{Bo86,Bo88} on geodesic
currents on surfaces is that the notion of the geometric intersection
number between (free homotopy classes of) closed curves on a closed
hyperbolic surface $S_g$ extends to a continuous symmetric bilinear
map $i: Curr(S_g)\times Curr(S_g)\to\mathbb R$. A crucial feature of
this construction is that if $\eta_c$ is the current determined by a
free homotopy class $c$ of closed curves and if $L_\rho$ is the
Liouville currents corresponding to a hyperbolic structure $\rho$ on
$S_g$ then $i(L_\rho, \eta_c)=\ell_\rho(c)$ where $\ell_\rho$ is the
marked length spectrum corresponding to $\rho$. That is $i(L_\rho,
\eta_c)$ is the $\rho$-length of the shortest with respect to $\rho$
curve in the class $c$.

It turns out that in the context of a free group $F$ of finite rank
$k\ge 2$ there exists a natural ``intersection form''
\[
I: FLen(F)\times Curr(F)\to \mathbb R
\]
where $FLen(F)$ is the space of all hyperbolic length functions on $F$
corresponding to free and discrete isometric actions of $F$ on
$\mathbb R$-trees. The form $I$ is continuous, linear with respect to
the second argument and $\mathbb R$-homogeneous with respect to the
first argument. Also, the form $I$ is equivariant with respect to the
left action of $Out(F)$ on $FLen(F)$ and $Curr(F)$. Moreover, as in
the surface case, if $\eta_g$ is the current determined by the
conjugacy class $[g]$ of $g\in F$ and if $\ell\in FLen(F)$ is
arbitrary then $I(\ell, \eta_g)=\ell(g)$. A normalized version of the
intersection form $I$, as explained in Section~\ref{sect:distortion}
below, already appears in \cite{Ka} where it serves as the main tool
for computing the conjugacy distortion spectrum of free groups
automorphisms. The definition given in the present paper is due to
Lustig and Hubert~\cite{Lu}.

It turns out that some symmetries and dualities applicable to
hyperbolic surfaces break down for the case of currents of free
groups.  Thus in Theorem~\ref{asym} we prove that there does not exist
a natural symmetric extension of $I$ to a map $Curr(F)\times
Curr(F)\to \mathbb R$.  To do that we interpret the value $I(\ell,
n_A)$, where $\ell\in FLen(F)$ is arbitrary and where $n_A$ is the
uniform current on $F$ corresponding to a free basis $A$, as the
``generic stretching factor" $\lambda_A(\ell)$ of $\ell$ with respect
to $A$.  Here $\lambda_A(\ell)$ is approximated by the
$\ell$-distortion of a ``random" $A$-geodesic. That is, for a long
random cyclically reduced word $w$ of length $n$ over $A$ we have
$\ell(w)/n\approx \lambda_A(\ell)$. The obstruction to a symmetric
extension of $I$ is caused by the fact that there exists $\phi\in
Out(F)$ such that $\lambda_A(\phi \ell_A)\ne
\lambda_A(\phi^{-1}\ell_A)$ where $\ell_A$ is the length function on
$F$ corresponding to the action of $F$ on its Cayley graph with
respect to $A$.

There is one obvious and notable exception in the topics covered in
that we do not discuss the Patterson-Sullivan-Bowen-Margulis embedding
from the outer space into the space of geodesic currents on a free
group (see~\cite{Fur} for an excellent discussion Patterson-Sullivan
measures in the general word-hyperbolic group context). This topic is
addressed in detail in a subsequent paper of the author and Tatiana
Nagnibeda~\cite{KN}.

Because this paper is, in part, expository, some of the proofs are
omitted or just sketched and there is an emphasis on aspects that are
not addressed in detail in ~\cite{Ka}.

The author is very grateful to Vadim Kaimanovich, Martin Lustig,
Tatiana Smirnova-Ngnibeda, Jerome Los and Paul Schupp for the many
helpful discussions.  The author expresses special thanks to Centre de
Recerca Matem\`atica at Barcelona where most of this paper was
written, and to the organizers of the special year in Geometric Group
Theory at CRM, Enric Ventura and Jose Burillo, for the invitation to
visit.

\section{Basic definitions}
\begin{conv}
  For the remainder of the paper, unless specified otherwise, let $F$
  be a finitely generated free group of rank $k\ge 2$. We will denote
  by $\partial F$ the hyperbolic boundary of $F$ in the sense of the
  theory of word-hyperbolic groups. Since $F$ is free, $\partial F$
  can also be viewed as the space of ends of $F$ with the standard
  ends-space topology.

  Thus $\partial F$ is a topological space homeomorphic to the Cantor
  set. We will also denote
\[
\partial^2 F:=\{(\zeta,\xi): \zeta,\xi\in \partial F \text{ and }
\zeta\ne \xi\}.
\]
\end{conv}

\begin{defn}[Geodesic Currents]
  Let $F$ be a free group of finite rank $k\ge 2$. A \emph{geodesic
    current} on $F$ is a positive locally finite Borel measure on
  $\partial^2 F$ that is $F$-invariant. We denote the space of all
  geodesic currents on $F$ by $Curr(F)$.
  
  The space $Curr(F)$ comes equipped with a weak topology: for $\nu_n,
  \nu\in Curr(F)$ we have $\displaystyle\lim_{n\to\infty}\nu_n=\nu$
  iff for every two disjoint closed-open sets $S,S'\subseteq \partial
  F$ we have $\displaystyle\lim_{n\to\infty}\nu_n(S\times
  S')=\nu(S\times S')$.
\end{defn}

\begin{defn}[Projectivized Geodesic Currents]
  For two nonzero geodesic currents $\nu_1, \nu_2\in Curr(F)$ we say
  that $\nu_1$ is equivalent to $\nu_2$, denoted $\nu_1\sim \nu_2$, if
  there exists a nonzero scalar $r\in \mathbb R$ such that $\nu_2=r
  \nu_1$. We denote \[\mathbb PCurr(F):=\{\nu\in Curr(F): \nu\ne
  0\}/\sim\] and call it the \emph{space of projectivized geodesic
    currents on $F$}. Elements of $\mathbb PCurr(F)$ (that is, scalar
  equivalence classes of elements of $Curr(F)$) are called
  \emph{projectivized geodesic currents}. The space $\mathbb PCurr(F)$
  is endowed with the quotient topology. We will denote the
  $\sim$-equivalence class of a nonzero geodesic current $\nu$ by
  $[\nu]$.
\end{defn}

We will see later on that, once a simplicial chart (defined below) on
$F$ is fixed, there is a natural embedding $\mathbb PCurr(F)\to
Curr(F)$ that provides a section to the natural quotient map
$Curr(F)-\{0\}\to \mathbb PCurr(F)$.

\begin{rem}[A note on the symmetrization]
  Denote by $\sigma:\partial^2 F\to \partial^2 F$ the \emph{flip} map
  $\sigma:(\zeta,\xi)\mapsto (\xi,\zeta)$ for $(\zeta,\xi)\in
  \partial^2 F$.
  
  Let $Curr_s(F)$ be the set of all $\nu\in Curr(F)$ such that $\nu$
  is $\sigma$-invariant. It is clear that $Curr_s(F)$ is a closed
  linear subspace of $Curr(F)$. The image of $Curr_s(F)$ in $\mathbb
  PCurr(F)$ is denoted by $\mathbb PCurr_s(F)$.
  
  Frequently the requirement for $\nu$ to be flip-invariant is
  included in the definition of a (projectivized) geodesic current. We
  will not do that here since most arguments and statements (at least
  those discussed in this paper) work in exactly the same way in both
  contexts and for the reasons of savings space and simplifying
  notations we prefer not to impose an extra condition at the
  definition level.
  
  We should also add that there is a natural retraction $r: Curr(F)\to
  Curr_s(F)$ defined as follows. If $\nu\in Curr(F)$ and $S\subseteq
  \partial^2 F$, we have $(r(\nu))(S)=\frac{1}{2}(\nu(S)+\nu(\sigma
  S))$.
\end{rem}

An important basic fact about geodesic currents is:
\begin{prop}
  The space $\mathbb PCurr(F)$ is compact.
\end{prop}

\begin{conv}
  If $\gamma$ is an edge-path or a circuit in some graph, we will
  denote by $|\gamma|$ the edge-length of $\gamma$. Similarly, if $w$
  is a word in some alphabet or a \emph{cyclic word} in some alphabet
  (to be defined later), we denote by $|w|$ the length of $w$, that
  is, the number of letters in $w$. If $\Delta$ is a graph, we will
  denote by $\mathcal P(\Delta)$ the set of all edge-paths of finite
  positive edge-length in $\Delta$.
\end{conv}

\begin{defn}[Simplicial charts]
  Let $\Gamma$ be a finite connected graph such that
  $\pi_1(\Gamma)\cong F$. Let $\alpha:F\to \pi_1(\Gamma,p)$ be an
  isomorphism, where $p$ is a vertex of $\Gamma$. We will call such
  $\alpha$ a \emph{simplicial chart} for $F$.
\end{defn}

\begin{conv}\label{conv:graph}
  Let $\alpha:F\to \pi_1(\Gamma,p)$ be a simplicial chart.
  
  Then $X:=\widetilde \Gamma$ is a topological tree. Denote the
  covering map from $X$ to $\Gamma$ by $q:X\to \Gamma$. For $\gamma\in
  \mathcal P(X)$ we will call the reduced edge-path $v=q(\gamma)$ in
  $\Gamma$ the \emph{label} of $\gamma$.
  
  Let $\partial X$ denote the space of ends of $X$ with the natural
  ends-space topology. Then we obtain a canonical $\alpha$-equivariant
  homeomorphism $\hat \alpha:\partial F\to\partial X$.
  
  This homeomorphism can be thought of in the following way. Suppose
  we endow $\Gamma$ with the structure of a metric graph, that is, we
  assign each edge of $\Gamma$ a positive length. This turns $X$ into
  an $\mathbb R$-tree with a discrete isometric action of
  $\pi_1(\Gamma,p)$.  Moreover, $X$ is quasi-isometric to $F$ and, if
  $F$ is equipped with a word metric and $p'$ is a lift of $p$ to
  $X=\widetilde G$ then the orbit map $\tilde\alpha:F\to X$,
  $\tilde\alpha:f\to \alpha(f)p'$, is a quasi-isometry. This
  quasi-isometry extends to a homeomorphism $\partial F\to \partial X$
  that is equal to $\hat\alpha$.

  If $\alpha$ is fixed, we will usually suppress explicit mention of
  $\hat \alpha$ and also of the map $\alpha$ itself when talking about
  the action of $F$ on $X$ and on $\partial F$ arising from this
  situation. Thus we also have an identification $\hat \alpha:
  \partial^2F\to \partial^2 X$. A crucial feature of this construction
  is that $\hat\alpha$ does not depend on the choice of a metric graph
  structure on $\Gamma$.
\end{conv}

\begin{conv}
  Suppose $\Gamma, \alpha, p$ and $X$ are as in
  Convention~\ref{conv:graph}.
  
  Similarly to the case of $F$, we denote by $\partial^2 X$ the set of
  all pairs $(\zeta_1,\zeta_2)$ such that $\zeta_1, \zeta_2\in
  \partial X$ and $\zeta_1\ne \zeta_2$. For $(\zeta_1,\zeta_2)\in
  \partial^2 X$ we denote by $[\zeta_1,\zeta_2]$ the simplicial
  (non-parameterized) geodesic from $\zeta_1$ to $\zeta_2$ in $X$.
  Thus $[\zeta_1,\zeta_2]$ is a subgraph of $X$ isomorphic to the
  simplicial line, together with a choice of direction on that line.
\end{conv}

\begin{defn}
  For every (oriented) reduced edge-path $\gamma$ in $X$ of positive
  edge-length denote

\begin{gather*}
  Cyl_\alpha(\gamma):= \{(x,y)\in \partial^2 F: \gamma\subseteq
  [\hat\alpha(x),\hat\alpha(y)] \text{ in } X \\ \text{ and the
    orientations on $\gamma$ and on $[\hat\alpha(x),\hat\alpha(y)]$
    agree}\}
\end{gather*}
Note that by definition $Cyl_\alpha(\gamma)$ is a subset of
$\partial^2 F$ (rather than of $\partial^2 X$). This distinction
becomes important when $\alpha$ and $\Gamma$ are not fixed but come
from different points in the outer space. However, when $\alpha$ is
fixed, we will often suppress the subscript $X$ and denote
$Cyl(\gamma)=Cyl_\alpha(\gamma)$.
\end{defn}

The collection of all sets $Cyl(\gamma)$, where $\gamma$ varies over
$\mathcal P(X)$, gives a basis of closed-open sets for $\partial^2 F$.
Hence it is easy to see that:

\begin{lem}\label{lem:conv}
  For $\nu_n, \nu\in Curr(F)$
  $\displaystyle\lim_{n\to\infty}\nu_n=\nu$ iff
  $\displaystyle\lim_{n\to\infty}\nu_n(Cyl(\gamma))=\nu(Cyl(\gamma))$
  for every reduced edge-path $\gamma$ in $X$ of positive edge-length.
  Moreover, for $\nu,\nu'\in Curr(F)$ we have $\nu=\nu'$ iff
  $\nu(Cyl(\gamma))=\nu'(Cyl(\gamma))$ for every $\gamma\in \mathcal
  P(X)$.
\end{lem}

Note that for any $f\in F$ and $\gamma\in \mathcal P(X)$ we have $f
Cyl(\gamma)=Cyl(f \gamma)$.  Since geodesic currents are, by
definition, $F$-invariant, for a geodesic current $\nu$ and for
$\gamma\in \mathcal P(X)$ the value $\nu(Cyl(\gamma))$ only depends on
the label $q(\gamma)$ of $\gamma$.

\begin{defn}[Number of occurrences of a path in a current]
  Let $\alpha:F\to \pi_1(\Gamma,p)$ be a simplicial chart and let
  $X=\widetilde \Gamma$.  For a path $v\in \mathcal P(\Gamma)$ and for
  $\nu\in Curr(F)$ we denote $\langle
  v,\nu\rangle_\alpha:=\nu(Cyl(\gamma))$ where $\gamma$ is any lift of
  $v$ to $X$. We call $\langle v,\nu\rangle_\alpha$ the \emph{number
    of occurrences of $v$ in $\nu$}.
\end{defn}

In view of Lemma~\ref{lem:conv} one can view a simplicial chart
$\alpha:F\to \pi_1(\Gamma,p)$ as providing a ``coordinate system" on
$Curr(F)$. A current $\nu\in Curr(F)$ is uniquely defined by its
collection of ``coordinates" $(\langle v,\nu\rangle_\alpha)_{v\in
  \mathcal P(\Gamma)}$. Each path $v\in \mathcal P(\Gamma)$ can be
thought of as defining the ``coordinate function" $\langle v,
\cdot\rangle_\alpha$ on $Curr(F)$.

As usual we will often omit the subscript $\alpha$ if the chart
$\alpha$ is fixed.

The following lemma is a well-known fact about word-hyperbolic groups
(see, for example~\cite{GH}).
\begin{lem}\label{qi}
  Let $G$ be a group acting properly discontinuously and co-compactly
  by isometries on a geodesic hyperbolic metric space $X$ (hence $G$
  is word-hyperbolic). Let $\phi$ be an automorphism of $G$. Thus the
  action of $G$ on $X$ extends to an action of $G$ on $\partial X$ and
  $\phi$ determines a quasi-isometry of $X$ that extends to a
  homeomorphism $\tilde\phi: \partial X\to\partial X$.

  Then for any $g\in G$ and $\xi\in \partial X$ we have
  $\phi(g)\tilde\phi(\xi)=\tilde\phi(g\xi)$
\end{lem}

Lemma~\ref{qi} is a crucial ingredient in defining the action of
$Out(F)$ on $Curr(F)$.

\begin{lem}\label{lem:action}
  Let $\phi$ be an automorphism of $F$ and let $\nu\in Curr(F)$.
  Define a measure $\phi\nu$ on $\partial^2 F$ by setting
  $\phi\nu(S):=\nu(\phi^{-1}(S))$ for a Borel subset $S\subseteq
  \partial^2 F$.  Then $\phi\nu$ is a geodesic current on $F$.
\end{lem}
\begin{proof}
  Let $S\subseteq\partial^2 F$ and $f\in F$. We need to check that
  $\phi\nu(fS)=\phi\nu(S)$.
  
  By definition $\phi\nu(fS)=\nu(\phi^{-1}(fS))$. By Lemma~\ref{qi}
  $\phi^{-1}(fS)=\phi^{-1}(f)\phi^{-1}(S)$. Since $\nu$ is
  $F$-invariant, we have
\[
\nu(\phi^{-1}(fS))=\nu(\phi^{-1}(f)\phi^{-1}(S))=\nu(\phi^{-1}(S))=\phi\nu(S).
\]
Thus $\phi\nu(fS)=\phi\nu(S)$, as required.
\end{proof}

\begin{prop}
  The map $(\phi,\nu)\mapsto \phi\nu$, where $\phi\in Aut(F), \nu\in
  Curr(F)$, defines a left action of $Aut(F)$ on $Curr(F)$ by
  continuous linear transformations. Moreover, $Inn(F)$ is contained
  in the kernel of this action, which, therefore, factors to the
  action of $Out(F)$ on $Curr(F)$.
\end{prop}
\begin{proof}
  It is clear from the definition that for a fixed $\phi\in Aut(F)$
  the map $\nu\mapsto \phi\nu$ is linear. Let us check the continuity
  of this map.  Suppose that $\nu_n, \nu\in Curr(F)$ and
  $\displaystyle\lim_{n\to\infty}\nu_n=\nu$. We need to show that
  $\displaystyle\lim_{n\to\infty}\phi\nu_n=\phi\nu$. Consider an
  arbitrary closed-open set $S\subseteq\partial^2 F$. Since $\phi$
  defines a homeomorphism of $\partial^2 F$, the set $\phi^{-1}(S)$ is
  also closed-open. Hence $\displaystyle\lim_{n\to\infty}\nu_n=\nu$
  implies that
  $\displaystyle\lim_{n\to\infty}\nu_n(\phi^{-1}(S))=\nu(\phi^{-1}(S))$,
  that is $\displaystyle\lim_{n\to\infty}\phi\nu_n(S)=\phi\nu(S)$.
  Since $S$ was an arbitrary closed-open set, this implies that
  $\displaystyle\lim_{n\to\infty}\phi\nu_n=\phi\nu$.
  
  Let us now check that $(\phi,\nu)\mapsto \phi\nu$ gives a left
  action of $Aut(F)$ on $Curr(F)$.  Suppose $\phi, \psi\in Aut(F)$ and
  $\nu\in Curr(F)$. We need to check that
  $(\phi\psi)(\nu)=\phi(\psi\nu)$. Take an arbitrary closed-open
  $S\subseteq\partial^2 F$. We have
  $[(\phi\psi)\nu](S)=\nu((\phi\psi)^{-1}(S))=\nu(\psi^{-1}\phi^{-1}(S))$.
  On the other hand
  $[\phi(\psi\nu)](S)=(\psi\nu)(\phi^{-1}(S))=\nu(\psi^{-1}\phi^{-1}(S))$,
  as required.
  
  Finally, observe that inner automorphisms act trivially on
  $Curr(F)$. Let $f\in F$ and consider an automorphism $\tau_f$ of $F$
  defines as $\tau_f(g)=fgf^{-1}$ for $g\in F$. Note that for every
  point $\xi\in \partial F$ we have $\tau_f(\xi)=f\xi$ and
  $\tau_f^{-1}(\xi)=f^{-1}\xi$. Hence for any closed-open
  $S\subseteq\partial^2 F$ we have $\tau_f^{-1}(S)=f^{-1}(S)$. Thus
  for any geodesic current $\nu$ and any Borel set $S$ as above
\[
(\tau_f\nu)(S)=\nu(\tau_f^{-1}S)=\nu(f^{-1}S)=\nu(S)
\]
where the last equality holds by $F$-invariance of $\nu$. Thus
$\tau_f\nu=\nu$ and we see that inner automorphisms of $F$ act
trivially on $Curr(F)$, as required.
\end{proof}

\begin{rem}
  Suppose that $\alpha, \beta: F\to\pi_1(\Gamma,p)$ are simplicial
  charts such that $\alpha^{-1}\beta$ is an inner automorphism of $F$
  and let $X=\widetilde\Gamma$. Then the maps $\hat\alpha, \hat\beta:
  \partial^2F\to \partial^2 X$ differ by a translation by an element
  of $F$.  That is, there is $g\in F$ such that for every
  $(\zeta,\xi)\in \partial^2 F$ we have
  $(\hat\alpha^{-1}\hat\beta)(\zeta,\xi)=(g\zeta,g\xi)$. Thus we have
  $Cyl_\alpha(\gamma)=g Cyl_\beta(\gamma)$ for every $\gamma\in
  \mathcal P(X)$. Hence for every $\nu\in Curr(F)$ and every path
  $v\in \mathcal P(\Gamma)$ we have $\langle v;\nu\rangle_\alpha=
  \nu(Cyl_\alpha(v))=\nu(Cyl_\beta(v))=\langle v;\nu\rangle_\beta$.
\end{rem}

\section{Projectivization and related questions}

\begin{conv}
  We will denote by $FLen(F)$ the space of all hyperbolic length
  functions $\ell:F\to\mathbb R$ corresponding to free and discrete
  actions of $F$ on $\mathbb R$-trees. We will denote by $\mathbb
  PFLen(F)$ or by $CV(F)$ the space of projective equivalence classes
  of nonzero elements of $FLen(F)$. Here two functions in $FLen(F)$
  are equivalent if they are scalar multiples of each other. The space
  $FLen(F)$ comes equipped with the weak topology of pointwise
  convergence on finite subsets of $F$. The space $CV(F)$ inherits the
  quotient topology. We will denote the projective equivalence class
  of $\ell\in FLen(F)$ by $[\ell]$.
  
  It is well-known that $CV(F)$ is precisely the \emph{Culler-Vogtmann
    outer space} of $F$, as defined in \cite{CV}, and we shall exploit
  both points of view here. Thus an element $\ell'\in CV(F)$ can also
  be represented as a simplicial chart $\alpha: F\to \pi_1(\Gamma,p)$
  where $\Gamma$ has no degree-one or degree-two vertices and is
  equipped with the structure of a metric graph so that the sum of the
  length of unoriented edges of $\Gamma$ is equal to $1$. Then
  $X=\widetilde \Gamma$ is an $\mathbb R$-tree with a free and
  discrete isometric action of $F$ on $X$ via $\alpha$. and the Let
  $\ell$ be the hyperbolic length function for this action of $F$ on
  $X$. Then $[\ell]=\ell'$.
  
  If $A$ is a free basis of $F$, we denote by $\ell_A\in Len(F)$ the
  hyperbolic length function on $F$ corresponding to the action of $F$
  on its Cayley graph with respect to $A$. Thus for $g\in F$
  $\ell_A(g)$ is the cyclically reduced length of $g$ with respect to
  $A$. We will denote the freely reduced length of $g$ with respect to
  $A$ by $|g|_A$.
\end{conv}

\begin{conv}
  Suppose that $F$ acts minimally, freely and discretely on $\mathbb
  R$-tree $X$ and that $\ell:F\to \mathbb R$ is a hyperbolic length
  function associated to this action.
  
  Then $\Gamma=X/F$ is a finite metric graph. If $e$ is an edge of $X$
  or of $\Gamma$, we will still denote the length of $e$ by $\ell(e)$.
\end{conv}

\begin{rem}[Left and Right actions on the outer space]
  The group $Out(F)$ acts on both $FLen(F)$ and $CV(F)$ by
  homeomorphisms. The traditional action is the \emph{right} action
  that, at the level of length functions, is given as follows.
  
  Let $\ell\in FLen(F)$ and let $\phi\in Aut(F)$ be an automorphism
  representing its outer automorphism class $[\phi]\in Out(F)$. Then
  $\ell [\phi]:=(\ell\circ \phi): F\to\mathbb R$.
  
  However, the natural action of $Out(F)$ on $Curr(F)$ and $\mathbb
  PCurr(F)$, which we described above, is a \emph{left} action.
  Hence, for the equivariance and embeddability purposes, in this
  paper we consider the \emph{left} action of $Out(F)$ on $FLen(F)$
  and $CV(F)$. In the above notations it is defined as follows:
  $[\phi]\ell:=(\ell\circ \phi^{-1}): F\to\mathbb R$. Thus $[\phi]\ell
  (g)= \ell(\phi^{-1}g)$ for $g\in F$. It is easy to see that this
  indeed defines left actions of $Out(F)$ on $FLen(F)$ and $CV(F)$.
  
  The left action is ``natural" in the following sense. There is a
  natural left action of $Aut(F)$ on elements of $F$ and on free bases
  of $F$. Let $\phi\in Aut(F)$ and let $A$ be a free basis of $F$.
  Let $\ell_A:F\to \mathbb R$ be the hyperbolic length function
  corresponding to the left action of $F$ on its Cayley graph with
  respect to $A$. Thus $\ell_A(g)$ is the cyclically reduced length of
  $g$ over $A$. Then it is clear that
  $\ell_{\phi(A)}(g)=\ell_A(\phi^{-1}g)$. Thus under our convention
  $\phi \ell_A =\ell_{\phi(A)}$.
\end{rem}

\begin{defn}[Normalization]\label{defn:weight}
  Let $\alpha: F\to \pi_1(\Gamma,p)$ be a simplicial chart for $F$ and
  $X=\widetilde \Gamma$. Let $\nu\in Curr(F)$ be a geodesic current.
  Denote by $E\Gamma$ the set of all oriented edges of $\Gamma$. Put
\[
\omega_\alpha(\nu):=\sum_{e\in E\Gamma} \langle e, \nu\rangle_\alpha.
\]
We call $\omega_\alpha(\nu)$ the \emph{weight of $\nu$ with respect to
  $\alpha$}. For a nonzero $\nu\in Curr(F)$ denote $\nu_\alpha:=
\nu/\omega_\alpha(\nu)$. Thus $[\nu]=[\nu_\alpha]$ and $\nu_\alpha$ is
the unique scalar multiple of $\nu$ that has $\alpha$-weight $1$. We
call $\nu_\alpha$ the \emph{$\alpha$-normalized representative of
  $\nu$} and, in general, we will say that a current is
\emph{$\alpha$-normalized} if it has $\alpha$-weight $1$.
\end{defn}

The following lemma is an easy exercise. It gives an explicit
criterion for the convergence of projectivized currents.

\begin{lem}
  Let $\alpha: F\to \pi_1(\Gamma,p)$ be a simplicial chart for $F$.
  
  Then for any nonzero $\nu, \nu_n\in Curr(F)$ (where $n=1,2,\dots$)
  we have:\newline $\displaystyle\lim_{n\to\infty} [\nu_n]=[\nu]$ in
  $\mathbb PCurr(F)$ if and only if $\displaystyle\lim_{n\to\infty}
  (\nu_n)_\alpha=\nu_\alpha$ in $Curr(F)$.
  
  The map $i_\alpha: \mathbb PCurr(F)\to Curr(F)$, $[\nu]\mapsto
  \nu_\alpha$, is an $Out(F)$-equivariant topological embedding of
  $\mathbb PCurr(F)$ in $Curr(F)$.
\end{lem}

Let $\alpha: F\to \pi_1(\Gamma,p)$ be a simplicial chart for $F$ and
$X=\widetilde \Gamma$. Let $g\in F$ be a nontrivial element and let
$c=c(g)$ be the cyclic path in $\Gamma$ representing $[g]$. Note that

\[
|c|=\sum_{e\in E\Gamma} \langle e, c\rangle=\sum_{e\in E\Gamma}
\langle e, [g]\rangle=\sum_{e\in E\Gamma} \langle e,
\eta_g\rangle=\omega_\alpha(\eta_g)
\]

Thus for the $\alpha$-normalized rational current $(\eta_g)_\alpha$ we
have $(\eta_g)_\alpha=\frac{\eta_g}{|c|}$. Hence for every path $v\in
\mathcal P(\Gamma)$ we have
\[
\langle v,(\eta_g)_\alpha\rangle= \langle
v,\eta_g\rangle/\omega_\alpha(\eta_g)=\langle
v,(\eta_g)\rangle/|c(g)|.
\]

So $\langle v,(\eta_g)_\alpha\rangle$ is equal to the number of
occurrences of $v$ in $c(g)$ divided by the length of $c(g)$. This
motivates the following:

\begin{defn}[Frequencies]
  Let $\alpha: F\to \pi_1(\Gamma,p)$ be a simplicial chart for $F$.
  For an arbitrary nonzero current $\nu\in Curr(F)$ and for any path
  $v\in \mathcal P(\Gamma)$ we call $\langle v,(\nu)_\alpha\rangle=
  \langle v,\nu\rangle/\omega_\alpha(\nu)$ the \emph{frequency} of $v$
  in $\nu$.
\end{defn}

Similarly to the situation in $Curr(F)$, once a simplicial chart
$\alpha$ is fixed, it can be thought of as providing a ``coordinate
system" on $\mathbb PCurr(F)$. Every projective class $[\nu]\in
\mathbb PCurr(F)$ is uniquely determined by its ``frequency
coordinates" $(\langle v,(\nu)_\alpha\rangle)_{v\in \mathcal
  P(\Gamma)}$. This point of view is explored in detail in \cite{Ka}.

\section{Other models}

There are several other spaces that are naturally homeomorphic to
$Curr(F)$ and $\mathbb PCurr(F)$.

We will briefly discuss them here omitting most of the details.

\begin{defn}[Semi-infinite sequences and one-sided shift]
  Let $\alpha:F\to \pi_1(\Gamma,p)$ be a simplicial chart for $F$ and
  let $X=\widetilde \Gamma$. Let $\Omega(\Gamma)$ be the set of all
  semi-infinite reduced edge-paths in $\Gamma$. For each vertex $a\in
  V\Gamma$ denote by $\Omega_a(\Gamma)$ the set of all $\gamma\in
  \Omega(\Gamma)$ that start at $a$. Thus $\Omega(\Gamma)=\sqcup_{a\in
    V\Gamma}\Omega_a(\Gamma)$. Each $\Omega_a(\Gamma)$ is identified
  in the obvious way with $\partial X$ and topologized accordingly,
  making it into a Cantor set. We give $\Omega_a(\Gamma)$ the topology
  of the disjoint union of several Cantor sets.
  
  Let $T_\Gamma: \Omega(\Gamma)\to \Omega(\Gamma)$ be the \emph{shift
    map}, that erases the first edge of each $\gamma\in
  \Omega(\Gamma)$. Then $T_\Gamma$ is easily seen to be continuous.
  
  Let $\mathcal S(\Gamma)$ be the space of all positive Borel measures
  $\mu$ on $\Omega(\Gamma)$ that are $T_\Gamma$-invariant, that is,
  have the property that for every Borel set $A\subseteq
  \Omega(\Gamma)$ we have $\mu(A)=\mu(T_\Gamma^{-1}A)$. Let $\mathbb
  PS(\Gamma)$ be the set of all $\mu\in S(\Gamma)$ such that
  $\mu(\Omega(\Gamma))=1$, that is, $\mu$ is a probability measure.
\end{defn}

\begin{defn}[Bi-infinite sequences and two-sided shift]
  Let $\alpha:F\to \pi_1(\Gamma,p)$ be a simplicial chart for $F$ and
  let $X=\widetilde \Gamma$.  Denote by $\Sigma(\Gamma)$ the set of
  all maps $\varsigma: \mathbb Z\to E\Gamma$ such that for every $i\in
  \mathbb Z$ we have $\varsigma(i)\varsigma(i+1)\in \mathcal
  P(\Gamma)$, that is, $\varsigma(i)\varsigma(i+1)$ is a reduced
  edge-path in $\Gamma$. We give $\Sigma(\Gamma)$ the natural weak
  topology of pointwise convergence on all finite subintervals of
  $\mathbb Z$ which makes $\Sigma(\Gamma)$ compact. The space
  $\Sigma(\Gamma)$ comes equipped with a natural shift action of
  $\mathbb Z$: for each $n\in \mathbb Z$ and $\varsigma\in
  \Sigma(\Gamma)$ we have $(\tau_n\varsigma)(i):=\varsigma(i+n)$. Then
  $\tau_n$ is a homeomorphism of $\Sigma(\Gamma)$ and
  $\tau_n\tau_m=\tau_{n+m}$ for every $n,m\in \mathbb Z$.
  
  We denote by $\mathcal T(\Gamma)$ the space of all positive Borel
  measures on $\Sigma(\Gamma)$ that are invariant with respect to this
  shift action of $\mathbb Z$.
  
  If we choose a lift $\tilde V$ of the vertex set $V\Gamma$ to $X$,
  we can also think of $\Sigma(\Gamma)$ as the set of all bi-infinite
  geodesic paths in $X$ that at time $0$ pass through one of the
  elements of $\tilde V$

\end{defn}
\begin{defn}[Geodesic flow]
  Let $\alpha:F\to \pi_1(\Gamma,p)$ be a simplicial chart for $F$ and
  let $X=\widetilde \Gamma$. We now give both $\Gamma$ and $X$
  simplicial metrics with every edge having unit length, so that $X$
  becomes an $\mathbb R$-tree and $\Gamma$ becomes a metric graph.
  
  Define the \emph{geodesic flow space} $\mathcal G(\Gamma)$ as the
  set of all isometric embeddings $\gamma:\mathbb R\to X$. The set
  $\mathcal G(\Gamma)$ is endowed with the compact-open topology,
  which in this case coincides with the weak (pointwise convergence)
  topology.
  
  There is an obvious $F$-action and $\mathbb R$-action on $\mathcal
  G(\Gamma)$ by homeomorphisms defined as follows. For $g\in F$ and
  $\gamma:\mathbb R\to X$ in $\mathcal G(\Gamma)$ put $g\gamma:=g\circ
  \gamma$, so that for each $r\in \mathbb R$ we have $(g\gamma)(r)=g
  \gamma(r)$.  Similarly, if $t\in \mathbb R$ then $\varrho_t:\mathcal
  G(\Gamma)\to \mathcal G(\Gamma)$ is defined as
  $(\varrho_t\gamma)(r)=\gamma(r+t)$. Then $\varrho_0=Id$ and
  $\varrho_t\circ \varrho_s=\varrho_{t+s}$ for any $t,s\in \mathbb R$.
\end{defn}

Note that the quotient by the shift action $\mathcal G(\Gamma)/\mathbb
R$ is equal to $\partial^2 X$.

\begin{defn}[Bi-invariant measures on $\mathcal G(\Gamma)$]
  We denote by $BI(\Gamma)$ the space of all positive Borel measures
  on $\mathcal G(\Gamma)$ that are both $\mathbb R$- and
  $F$-invariant. This space comes equipped with the weak topology.

  We denote by $\mathbb P\mathcal{BI}(\Gamma)$ the space of projective
  equivalence classes of nonzero measures from $\mathcal{BI}(\Gamma)$,
  equipped with the quotient topology.
\end{defn}

\begin{prop}
  Let $\alpha:F\to \pi_1(\Gamma,p)$ be a simplicial chart for $F$ and
  let $X=\widetilde \Gamma$.  Then
  
  (a) The space $Curr(F)$ is naturally homeomorphic to $\mathcal
  S(\Gamma)$, to $\mathcal{BI}(\Gamma)$ and to $\mathcal T(\Gamma)$.
  
  (b) The space $\mathbb PCurr(F)$ is naturally homeomorphic to
  $\mathbb P\mathcal S(\Gamma)$, to $\mathbb P\mathcal{BI}(\Gamma)$
  and to $\mathbb P\mathcal T(\Gamma)$.
\end{prop}

\begin{proof}

  We will discuss briefly only the proof of (a) since part (b) is
  completely analogous.
  
  For a path $v\in\mathcal P(\Gamma)$ let $Cyl(v)$ denote the set of
  all paths $\gamma\in \Omega(\Gamma)$ that start with $v$. Then
  $Cyl(v)$ is an open-closed subset of $\Omega(\Gamma)$ and the sets
  $Cyl(v)$, where $v$ varies over all $\mathcal P(\Gamma)$ generates
  the Borel $\sigma$-algebra for $\Omega(\Gamma)$.
  
  Suppose now that $\mu\in S(\Gamma)$ is a shift-invariant measure on
  $\Omega(\Delta)$. We define a measure $\hat \mu$ on $\partial^2 X$
  as follows. For each path $\beta\in\mathcal P(X)$ put
  $\hat\mu(Cyl(\beta)):=\mu(Cyl(v))$ where $v$ is the label of
  $\beta$. It is not hard to see that $\hat\mu$ is a measure on
  $\partial^2 X$ that, by construction, is $F$-invariant. The map
  $\mu\mapsto\hat\mu$ provides a homeomorphism $S(\Gamma)\to Curr(F)$
  that factors to a homeomorphism between the projectivized versions
  of these spaces.
  
  Consider now the geodesic flow space $\mathcal G(\Gamma)$. Choose
  $t\in \mathbb R$ and a path $\beta\in\mathcal P(X)$. Define
  $Cyl(\beta,t)$ to be the set of all geodesics $\gamma:\mathbb R\to
  X$ such that $\gamma([t,t+|\beta|])=\beta$ and $\gamma$ maps
  $[t,t+|\beta|]$ to $\beta$ preserving the orientation. Then
  $Cyl(\beta,t)\subset \mathcal G(\Gamma)$ is closed-open and compact.
  The sets $Cyl(\beta,t)$, where $t\in\mathbb R, \beta\in\mathcal
  P(X)$, generate the Borel $\sigma$-algebra for $\mathcal G(\Gamma)$.
  
  Suppose now that $\nu\in Curr(F)$ is a geodesic current. We define a
  measure $\tilde\nu$ on $\mathcal G(\Gamma)$ as follows. For every
  cylinder $Cyl(\beta,t)$ put
  $\tilde\nu(Cyl(\beta,t)):=\nu(Cyl(\beta))$. Again, it is easy to see
  that $\tilde\nu$ is a measure which is both $F$- and $\mathbb
  R$-invariant, so that $\tilde\nu\in B(\Gamma)$. Also, the map
  $\nu\mapsto \tilde\nu$ is a homeomorphism between $Curr(F)$ and
  $B(\Gamma)$ that factors through to a homeomorphism between their
  projectivizations.

  Finally, let us discuss the identification between $Curr(F)$ and
  $\mathcal T(\Gamma)$. In the context of $\Sigma(\Gamma)$ the
  cylinder sets look as follows. Let $v\in \mathcal P(\gamma)$ be a
  path with $|v|=n$ and let $t\in \mathbb Z$.  Then $Cyl_t(v)$ is
  defined as the set of all $\varsigma\in Sigma(\Gamma)$ such that
  $\varsigma(i) \varsigma(i+1)\dots \varsigma(i+n-1)=v$. Again, the
  sets $Cyl_t(v)$, where $v\in \mathcal P(\Gamma)$ and $t\in \mathbb
  Z$, generate the Borel $\sigma$-algebra for $\Sigma(\Gamma)$.
  
  Suppose $\nu\in Curr(F)$. We associate to $\nu$ a Borel measure
  $\mu$ on $\Sigma(\Gamma)$ as follows. For any cylinder $Cyl_t(v)$
  set $\mu(Cyl_t(v)):=\langle v, \nu\rangle=\nu(Cyl(\gamma))$ where
  $\gamma$ is any lift of $v$ to $X$. It is not hard to see that $\mu$
  is shift-invariant, so that $\mu\in \mathcal T(\Gamma)$. This
  determines a map from $Curr(F)$ to $\mathcal T(\Gamma)$ that is
  easily seen to be bijective and continuous.
\end{proof}

The identification between $\mathbb PCurr(F)$ and the space $\mathbb
P\mathcal S(\Gamma)$ is particularly useful, for example, for proving
the density of rational currents in $Curr(F)$ and in $\mathbb
PCurr(F)$. Indeed, it is not hard to see that the elements of $\mathbb
P\mathcal T(\Gamma)$ corresponding to rational currents are precisely
the shift-invariant probability measures supported on the
$T_\Gamma$-periodic orbits in $\Omega(\Gamma)$. This connection is
explored in more detail in \cite{Ka}.

\section{Currents determined by conjugacy classes and the intersection form}

Let $g\in F$ be a nontrivial element. It canonically defines a pair of
distinct points $g^{\infty}, g^{-\infty} \in \partial F$, where
$g^{\infty}=\displaystyle\lim_{n\to\infty} g^n$ and
$g^{-\infty}=\displaystyle\lim_{n\to\infty} g^{-n}$. Note that
$(fgf^{-1})^{\infty}=fg^{\infty}$ and
$(fgf^{-1})^{-\infty}=fg^{-\infty}$ for every $f\in F$.

\begin{defn}[Rational currents]
  Let $g\in F$ be a nontrivial element that is not a proper power.
  Define a Borel measure $\eta_g$ on $\partial^2 F$ as follows.  For a
  closed-open subset $S\subseteq \partial^2 F$ let $\eta_g(S)$ be the
  number of those $F$-translates of $(g^{\infty},g^{-\infty})$ that
  belong to $S$. Obviously, $\eta_g$ is $F$-invariant, so that
  $\eta_g\in Curr(F)$.
  
  In view of the remark above $\eta_g(S)$ is equal to the number of
  points of the form $f(g^{\infty},g^{-\infty})f^{-1}$, where $f\in
  F$, that belong to $S$. Thus $\eta_g$ only depends on the conjugacy
  class of $g$ in $F$.
  
  Suppose now that $g$ is an arbitrary nontrivial element of $F$.
  Then we can uniquely represent $g$ as $g=h^t$ where $t\ge 1$ is an
  integer and $h\in F$ is not a proper power. We define $\eta_g:=t
  \eta_h$. Again, we see that $\eta_g$ only depends on the conjugacy
  class of $g$.
  
  Finally, for any nontrivial conjugacy class $[g]$ in $F$ define
  $\eta_{[g]}:=\eta_g$, where $g\in [g]$ is an arbitrary element.
  
  We shall refer to multiplies of currents of the form $\eta_{[g]}$,
  where $[g]$ is a nontrivial conjugacy class in $F$, as
  \emph{rational} currents.
  
  For any nontrivial $g\in F$ we denote the projective class
  $[\eta_g]$ of $\eta_g$ in $\mathbb PCurr(F)$ by $\mu_g$ or by
  $\mu_{[g]}$. Note that by definition if $n\ge 1$ then
  $\mu_g=\mu_{g^n}$.
\end{defn}

Suppose now that $\partial F$ is identified with $\partial
X=\partial\widetilde \Gamma$, as in the previous section. Let $w$ be a
nontrivial root-free conjugacy class in $F$. Then $w$ is represented
by a unique cyclically reduced closed circuit $w'$ in $\Gamma$. Choose
a particular cyclically reduced path $\gamma$ in $\Gamma$ representing
$w'$ (which can be done by choosing a vertex on $w'$).

Then it is easy to see the for every path $\gamma\in \mathcal P(X)$
$\eta_g(Cyl(\gamma))$ is equal to the number of those $(x,y)\in
\partial^2 X$ such that the geodesic $[x,y]$ contains $\gamma$ (with
the agreement of orientations) and such that $[x,y]$ is labelled by a
bi-infinite power of $\gamma$:
\[
\dots \gamma\gamma\gamma \dots
\]

\begin{defn}[Cyclic paths and cyclic words]
  A \emph{cyclic path} or \emph{circuit} in $\Gamma$ is an immersion
  graph-map $c:\mathbb S\to \Gamma$ from a simplicially subdivided
  oriented circle $\mathbb S$ to $\Gamma$. Let $u$ be an edge-path in
  $\Gamma$. An \emph{occurrence of $u$ in $c$} is a vertex of $\mathbb
  S$ such that, going from this vertex in the positive direction along
  $\mathbb S$, there exists an edge-path in $\mathbb S$ (not
  necessarily simple and not necessarily closed) which is labelled by
  $u$, that is, which is mapped to $u$ by $c$. We denote by $\langle
  u,c\rangle$ the number of occurrences of $u$ in $c$.
  
  If $A$ is a free basis of $F$ and $\Gamma$ is a bouquet of edges
  labelled by the elements of $A$, then a cyclic path in $\Gamma$ can
  also be thought of as a \emph{cyclic word} over $A$. A \emph{cyclic
    word} is an equivalence class of cyclically reduced words, where
  two cyclically reduced words are equivalent if they are cyclic
  permutations of each other.
\end{defn}

\begin{notation}
  Let $\alpha:F\to \pi_1(\Gamma,p)$ be a simplicial chart for $F$ and
  let $X=\widetilde \Gamma$. It is easy to see that every nontrivial
  conjugacy class $[g]$ in $F$ is uniquely represented by a reduced
  circuit in $\Gamma$ which in turn is uniquely represented by a
  cyclic path $c=c_\alpha(g)$ in $\Gamma$.
\end{notation}

The following is an immediate corollary of the definitions:

\begin{lem}\label{count}
  For every $\gamma\in \mathcal P(X)$ with label $v=q(\gamma)$ and for
  every nontrivial conjugacy class $[g]$ in $F$ with $c=c(g)$ we have
\[
\eta_{[g]}(Cyl(\gamma))=\langle v,c\rangle,
\]
that is,
\[
\langle v, \eta_{[g]}\rangle=\langle v,c\rangle
\]
\end{lem}

A similar statement holds at the level of frequencies. For example,
suppose that $(\alpha,\Gamma)$ is the bouquet of circles corresponding
to a free basis $A$ of $F$, and that $w$ is a cyclic word over $A$
representing $g$. If $v$ is a freely reduced word over $A$ then the
frequency of $v$ in $\eta_g$ is equal to the number of occurrences of
$v$ in $w$ divided by the length of $w$, which is the frequency of $v$
in $w$.

\begin{prop}\label{good}
  Let $g\in F$ be a nontrivial element and let $\phi\in Aut(F)$. Then
  $\phi\eta_g=\eta_{\phi(g)}$.
\end{prop}

\begin{proof}
  Let $S\subseteq \partial^2 F$ be a closed-open set and let
  $S'=\phi^{-1}(S)$. Recall that for $h\ne 1$ the measure $\eta_h$
  counts the number of points that are $F$-translates (or
  $F$-conjugates) of the pair $(h^{-\infty},h^{\infty})$ in a set.
  
  Thus $\phi\eta_g (S):=\eta_g(\phi^{-1}(S))$ is the number of
  $F$-translates of $(g^{-\infty},g^{\infty})$ in $S'=\phi^{-1}(S)$.
  Similarly, $\eta_{\phi(g)}(S)$ is the number of $F$-translates of
  the point $(\phi(g)^{-\infty},\phi(g)^{\infty})$ in $S$.
  Lemma~\ref{qi} implies that $\phi$ maps bijectively $F$-translates
  of $(g^{-\infty},g^{\infty})$ in $\phi^{-1}(S)$ to $F$-translates of
  $(\phi(g)^{-\infty},\phi(g)^{\infty})$ in $S$. Hence $\phi\eta_g
  (S)=\eta_{\phi(g)}(S)$ and, since $S$ was arbitrary,
  $\phi\eta_g=\eta_{\phi(g)}$.
\end{proof}

\begin{notation}
  Let $C=C(F)$ denote the set of all nontrivial conjugacy classes in
  $F$ and let $C_0=C_0(F)$ the the set of all elements of $C$ that are
  not proper powers. Denote by $r:C\to Curr(F)$ the map $r:[g]\mapsto
  \eta_{[g]}$ and denote by $\hat r: C\to \mathbb PCurr(F)$ the map
  $\hat r:[g]\mapsto \mu_{[g]}$.
\end{notation}

We summarize some of the properties of $r$ and $\hat r$ in the
following proposition (see \cite{Ka} for detailed arguments):

\begin{prop} We have:
\begin{enumerate}
\item The set $\hat r(C)=\hat r(C_0)$ is dense in $\mathbb PCurr(F)$
  and the set of all scalar multiples of elements of $r(C)$ is dense
  in $Curr(F)$.
  
\item The map $r:C\to Curr(F)$ is injective and the map $\hat
  r|_{C_0}:C_0\to \mathbb PCurr(F)$ is injective; consequently, the
  actions of $Out(F)$ on $Curr(F)$ and $\mathbb PCurr(F)$ are
  effective.
\end{enumerate}
\end{prop}

We can now define a natural ``intersection form" $I: FLen(F)\times
Curr(F)\to \mathbb R$ that will be $Out(F)$-equivariant and linear in
the second coordinate. Moreover, it will be ``natural" in the sense
that $I(\eta_g,l)=\ell(g)$ for every $g\in F$.

\begin{defn}[Intersection form]
  
  Let $\ell\in FLen(F)$ be a length function and let $\nu\in Curr(F)$
  is a geodesic current. Let $\alpha: F\to \pi_1(\Gamma,p)$ be a
  simplicial chart and let $\Gamma$ be given a metric graph structure
  so that $\ell$ is the hyperbolic length function corresponding to
  the action of $F$ on the $\mathbb R$-tree $X=\widetilde \Gamma$ via
  $\alpha$, where the metric on $X$ correspnds to the metric graph
  structure on $\Gamma$.  We define
\[
I(\ell,\nu):=\sum_{e\in E\Gamma} \ell(e)\langle e, \nu\rangle_\alpha,
\]
where $E\Gamma$ is the set of oriented edges of $\Gamma$.
\end{defn}

\begin{prop}
  The map $I: FLen(F)\times Curr(F)\to \mathbb R$ is continuous,
  linear with respect to the second argument and $Out(F)$-equivariant
  with respect to the left diagonal action of $Out(F)$ on
  $FLen(F)\times Curr(F)$. Also, $I$ is homogeneous with respect to
  the first argument, that is $I(r\ell,\nu)=r I(\ell,\nu)$ for every
  $r\ge 0$ and every $\ell\in FLen(F), \nu\in Curr(F)$.
  
  Moreover, for every $g\in F$, and $\ell\in FLen(F)$ we have
  $I(\ell,\eta_g)=\ell(g)$.

\end{prop}
\begin{proof}
  It is clear from the definition that $I: FLen(F)\times Curr(F)\to
  \mathbb R$ is linear with respect to the second argument,
  homogeneous with respect to the first argument, and, moreover, for
  any $\ell\in FLen(F)$ the function $I(\ell, -)$ is continuous on
  $Curr(F)$.
  
  Suppose that $\ell\in FLen(F)$ is a length function and that $\nu\in
  Curr(F)$ is a geodesic current. The length function $\ell$ defines a
  free and discrete isometric action of $F$ on an $\mathbb R$-tree $X$
  such that $\Gamma=X/F$ is a finite metric graph and that
  $X=\widetilde \Gamma$. Let $E$ be the collection of oriented edges
  of $\Gamma$. Let $g\in F$ be a nontrivial element. Let $c=c(g)$ be
  the reduced cyclic path in $\Gamma$ representing $[g]$. Then
  $\ell(g)$ is equal to the sum of the length of the edges of $c$,
  that is $\ell(g)=\sum_{e\in E} \ell(e)n_{c}(e)$. By
  Lemma~\ref{count} we have $n_{c}(e)=n_{\eta_g}(e)$. Thus
\[
\ell(g)=\sum_{e\in E} \ell(e)n_{c}(e)=\sum_{e\in E}
\ell(e)n_{\eta_g}(e)=I(\ell,\eta_g),
\]
as claimed.

Assuming the global continuity of $I$, let us check its
$Out(F)$-equivariance. Suppose $\ell\in FLen(F)$ and $g\in F$ and let
$\phi\in Aut(F)$. Recall that $\phi\eta_g=\eta_{\phi(g)}$.  Also, by
definition of the left action of $Out(F)$ on $FLen(F)$ we have $\phi
l= \ell\circ\phi^{-1}$. Thus we have
\[
I(\phi \ell, \phi \eta_g)=I(\phi \ell, \eta_{\phi(g)})=(\phi
\ell)(\phi(g))=(\ell\circ \phi^{-1})(\phi(g))=\ell(g)=I(\ell, \eta_g).
\]
Since the scalar multiples of rational currents are dense in
$Curr(F)$, the continuity of $I$ implies that $I(\phi l, \phi
\nu)=I(\ell,\nu)$ for every $\ell\in FLen(F)$ and every $\nu\in
Curr(F)$.

It remains to show that $I$ is continuous. We will give a sketch of
the argument here and leave some of the details to the reader. Let
$\ell\in CV, \nu\in Curr(F)$. We need to establish that $I$ is
continuous at $(\ell,\nu)$. Because $FLen(F)\times Curr(F)$ is
metrizable and locally compact, it suffices to prove that for every
two sequences $\ell_n\in FLen(F)$ and $\nu_n\in Curr(F)$ with
$\displaystyle\lim_{n\to\infty} \ell_n=l$ and
$\displaystyle\lim_{n\to\infty} \nu_n=\nu$ we have
$\displaystyle\lim_{n\to\infty} I(\ell_n,\nu_n)=I(\ell,\nu)$.

The length function $\ell$ defines a minimal action of $F$ on an
$\mathbb R$-tree $X$ with a finite quotient graph $\Gamma=X/F$.
Moreover, this action determines an isomorphism $\alpha:
F\to\pi_1(G,p')$ where $p'$ is a vertex of $\Gamma$.

Suppose first that $\Gamma$ has maximal possible number of edges among
all finite connected graphs with no degree-one and degree-two vertices
whose fundamental group is isomorphic to $F$. Then for all
length-functions $\ell'$ sufficiently close to $\ell$ we have
$\Gamma'=X'/F=X/F=\Gamma$ and $\alpha=\alpha'$ where $X'$ is the tree
corresponding to $\ell'$ and where $\alpha': F\to \pi_1(G,p')$ is
defined similarly to $\alpha$. Then the continuity of $I$ at
$(\ell,\nu)$ follows directly from the definition of $I$.

Suppose now that the number of edges of $\Gamma$ is not the maximal
possible. Then it suffices to consider the situation when the sequence
$\ell_n$ approximating $\ell$ has the following form.  There is a
finite graph $\Delta$ homotopy equivalent to $\Gamma$ and such that
$\Gamma$ is obtained from $\Delta$ by contracting to points a certain
(possibly empty) collection of edges $E'$ of $\Delta$. We will denote
this contraction by $\kappa:\Delta\to \Gamma$. There is an isomorphism
$\beta: F\to \pi_1(\Delta, p'')$, where $p''$ is a vertex of $\Delta$
such that $\kappa(p'')=p'$ and such that $\beta$ factors through
$\kappa_{\#}$ to $\alpha$, that is $\alpha=\kappa_{\#}\circ \beta$.
Each $\ell_n$ corresponds to making $\Delta$ into a finite metric
graph. Moreover, the original length function $\ell$ corresponds to a
semi-metric structure on $\Delta$, where every edge of $\Delta$ is
assigned a certain nonnegative length, with edges of $E'$ being given
zero length and the edges in $E\Delta-E'$ assigned the same length as
their images in $\Gamma$. Then again every conjugacy class $[g]$ in
$F$ is represented, via $\beta$, by a unique reduced cyclic path and
$\ell(g)$ is the length of that path. Note that $Y=\widetilde\Delta$
is no longer equivariantly homeomorphic to $X=\widetilde \Gamma$.
However, $X$ is obtained from $Y$ by performing an equivariant
collection of edge-contractions of the lifts of the edges of $E'$.
Moreover, for every $e\in E'$ we have $\ell_n(e)\to 0$ as
$n\to\infty$. Note that since $\nu_n$ converges to $\nu$, we have
$\displaystyle\lim_{n\to\infty}\nu_n(Cyl_\beta(e))=\nu(Cyl_\beta(e))$
for every edge $e$ of $\Delta$. In particular, for every such $e$ the
sequence $\nu_n(Cyl_\beta(e))$ is bounded. In case $e\in E'$ this
implies that $\displaystyle\lim_{n\to\infty}
\ell_n(e)\nu_n(Cyl_\beta(e))=0$.

The crucial point is that for every edge $e$ of $\Gamma$ (which we
still denote by $e$ when thought of as an edge of $\Delta$) we have
$Cyl_\alpha(e)=Cyl_\beta(e)\subseteq \partial^2 F$.  Consequently, for
every edge $e$ of $\Gamma$ and for any current $\mu\in Curr(F)$ we
have $\langle e, \mu\rangle_\beta=\langle e, \mu\rangle_\alpha$.

Hence
\begin{gather*}
  \lim_{n\to\infty}I(\ell_n,\nu_n)=\lim_{n\to\infty}\sum_{e\in
    E\Delta} \ell_n(e) \langle e,
  \nu_n\rangle_\beta=\\
  \lim_{n\to\infty}[\sum_{e\in E'} \ell_n(e) \langle e,
  \nu_n\rangle_\beta+\sum_{e\in EG} \ell_n(e) \langle e,
  \nu_n\rangle_\beta]=\\
  0+\lim_{n\to\infty}\sum_{e\in EG} \ell_n(e) \langle e, \nu_n\rangle_\beta=\\
  \lim_{n\to\infty}\sum_{e\in EG} \ell(e) \langle e,
  \nu\rangle_\alpha=I(\ell,\nu)
\end{gather*}
as required.
\end{proof}

\begin{rem}[Weight and length]
  The notion of weight is also natural in the following sense. Let
  $\alpha: F\to \pi_1(\Gamma,p)$ be a simplicial chart for $F$ and
  $X=\widetilde \Gamma$. Endow $X$ and $\Gamma$ with simplicial
  metrics and let $\ell_\alpha\in FLen(F)$ be the corresponding length
  function on $F$. In view of the definition of the intersection form
  $I(-,-)$ we see that for any nonzero current $\nu$ we have
  $I(\ell_\alpha, \nu_\alpha)=1$ and that $\nu_\alpha$ is the only
  representative of $[\nu]$ with this property.
\end{rem}

\section{Local formulas}

It is a natural and important question to understand what the
transition functions between coordinate systems on $Curr(F)$
corresponding to two different simplicial charts look like.

\begin{prop}\label{local}
  Let $\alpha: F\to \pi_1(\Gamma,p)$ and $\beta: F\to
  \pi_1(\Upsilon,s)$ be two simplicial charts for $F$. Let
  $X=\widetilde \Gamma$ and $Y=\widetilde \Upsilon$ be the
  corresponding topological trees.  There exists a constant
  $K=K(\alpha,\beta)>0$ with the following property.
  
  For any path $v\in \mathcal P(\Upsilon)$ there exist integers
  $C(u,v)=C(u,v,\alpha,\beta)\ge 0$ such that for every nontrivial
  $g\in F$ we have:
\[
\langle v, c_\beta(g)\rangle_\beta = \sum_{u\in \mathcal P(\Gamma),
  |u|\le K|v|} C(u,v) \langle u, c_\alpha(g)\rangle_\alpha.
\tag{$\ddag$}
\]

Therefore, since rational currents are dense in $Curr(F)$, for every
$\nu\in Curr(F)$ we have:
  \[
  \langle v, \nu\rangle_\beta = \sum_{u\in \mathcal P(\Gamma), |u|\le
    K|v|} C(u,v) \langle u, \nu\rangle_\alpha.
  \]
\end{prop}

\begin{proof}
  
  Since rational currents are dense in $Curr(F)$, by continuity it
  suffices to establish $(\ddag)$.
  
  This statement is proved in \cite{Ka} for the case where both
  $\Gamma$ and $\Upsilon$ are bouquets of $k$ loop-edges. We will
  refer to this as the bouquet-bouquet case. Note that this is the
  crucial case where the heart of the argument lies. The idea there is
  that the bouquet-bouquet case corresponds to considering a single
  automorphism $\phi$ of $F$ with respect to a fixed free basis $A$ of
  $F$. The statement of the proposition then says that for a fixed
  freely reduced word $v$ over $A$ and for any reduced cyclic word $w$
  over $A$ the number of occurrences of $v$ in $\phi(w)$ is an integer
  linear combination of the numbers of occurrences in $w$ of words $u$
  from a finite collection, where this collection depends only on
  $\phi$ and $v$ but not on $w$. This, in turn, can be easily
  established by induction on the word-length of an automorphism
  $\phi$, once the statement has been directly verified for the
  Nielsen automorphisms which generate $Aut(F)$.
  
  Since the bouquet-bouquet case is already covered, by a composition
  argument it suffices to prove statement for the case where $\Gamma$
  is a bouquet of $k$ loop-edges and $\Upsilon$ is arbitrary and for
  the case where $\Upsilon$ is a bouquet of $k$ loop-edges and
  $\Gamma$ is arbitrary.

  Consider first the case where $\Gamma$ is a bouquet of $k$ loop
  edges and $\Upsilon$ is arbitrary.  This case is in fact considered
  in \cite{Ka}, but we will repeat the argument.
  
  Choose a maximal tree $T$ in $\Upsilon$. Choose an orientation
  $E\Upsilon=E^+\sqcup E^-$ on $\Upsilon$. This defines a
  \emph{geometric basis $A_T$} of $F$ as follows. For each edge $e\in
  E^+-ET$ put $\gamma_e=[s,o(e)]_Te[t(e),s]\in \pi_1(\Upsilon,s)$.
  Note that $\gamma_e$ is a reduced edge-path from $s$ to $s$ in
  $\Upsilon$. Put $a_e=\beta^{-1}(\gamma_e)$ and put $A_T:=\{a_e| e\in
  E^+-ET\}$. Thus indeed $A_T$ is a free basis of $F$.
  
  By the bouquet-bouquet case we may assume that $(\Gamma, \alpha)$ is
  the the bouquet of edges corresponding to the free basis $A_T$ of
  $F$.
  
  There is an explicit way of rewriting a reduced cyclic word $w$ over
  $A_T$ into a reduced cyclic path $c_\Upsilon{w}$ as follows. Replace
  each $a_e^{\pm 1}$ in $w$ by $e^{\pm 1}$ and then between each pair
  of the sort $e_1^{\varepsilon} e_2^{\delta}$, where
  $\varepsilon,\delta\in\{\pm 1\}$, insert
  $[t(e_1^{\varepsilon}),o(e_2^{\delta})]$. The result is precisely
  the cyclic path $c_\Upsilon{w}$. It is now clear that every
  occurrence of a fixed path $v$ in $c_\Upsilon{w}$ must come from an
  occurrence of one of a finite collection of paths $u$ in $\Gamma$,
  where this collection depends only on $v$ but not on $w$.
  
  Note also that there is also an obvious way of rewriting a cyclic
  path $c$ in $\Upsilon$ into a cyclic word $w$ over $A_T$. Namely, we
  delete all the edges of $ET$ from $c$ and replace every $e^{\pm 1}$
  by $a_e^{\pm 1}$. This shows that the number of occurrences of a
  reduced word $z$ in $w$ is equal to the number of occurrences of the
  path $z'$ in $c$, where, again, $z'$ is obtained from $z$ by
  replacing each $a_e^{\pm 1}$ in $w$ by $e^{\pm 1}$ and then
  inserting between each pair of the sort $e_1^{\varepsilon}
  e_2^{\delta}$ the path $[t(e_1^{\varepsilon}),o(e_2^{\delta})]$.

  This shows, in addition, that the statement of the proposition also
  holds when $\Gamma$ and $\Upsilon$ as above exchange places. By the
  bouquet-bouquet case this implies that the statement holds when
  $\Gamma$ is arbitrary and when $\Upsilon$ is a bouquet of edges.
  
  This completes the proof of $(\ddag)$ and of the proposition.
\end{proof}

Note that in view of the discussion in Section~\ref{sect:approx} below
the condition $|u|\le K$ in the Proposition~\ref{local} can be
replaced by $|u|=K$.

\begin{thm}\label{thm:local}
  Let $\alpha: F\to \pi_1(\Gamma,p)$ be a simplicial chart for $F$ and
  let $\ell\in FLen(F)$ be a length function. Let $X=\widetilde
  \Gamma$.
  
  Then there exist an integer $K>0$ and some constants
  $d(u)=d(u,\alpha,\ell)\ge 0$, where $u$ varies over all reduced
  paths in $\Gamma$ of length at most $K$, with the following
  property.
  
  For any $\nu\in Curr(F)$ we have
\[
I(\ell,\nu)=\sum_{|u|\le K, u\in\mathcal P(\Gamma)} d(u)\langle u,
\nu\rangle_\alpha.
\]
In particular, for every nontrivial element $g\in F$
\[
I(\ell,\eta_g)=\ell(g)=\sum_{|u|\le K, u\in\mathcal P(\Gamma)}
d(u)\langle u, [g]\rangle_\alpha.
\]

\end{thm}

\begin{proof}
  Since the intersection form $I$ is continuous and rational currents
  are dense in $Curr(F)$, it suffices to establish the statement of
  the theorem for rational currents.

  The length function $\ell$ corresponds to a representation $\beta:
  F\to \pi_1(\Upsilon,s)$, where $\Upsilon$ is a metric graph and
  where $\ell$ is the hyperbolic length function associated to the
  action of $F$ (via $\beta$) on the $\mathbb R$-tree $Y=\widetilde
  \Upsilon$. We will still denote the length of an edge $e$ in
  $\Upsilon$ or in $Y$ by $\ell(e)$.
  
  Let $g\in F$ be an arbitrary nontrivial element. Let $c=c(g)$ be the
  reduced cyclic path in $\Upsilon$ representing $[g]$. Then
\[
\ell(g)=\sum_{e\in E\Upsilon} \ell(e) \langle e, c(g)\rangle.
\]

By Proposition~\ref{local} we have
\begin{gather*}
  I(\ell,\eta_g)=\ell(g)=\sum_{e\in E\Upsilon} \ell(e) \langle e,
  c(g)\rangle=\sum_{e\in E\Upsilon} \ell(e)\sum_{|u|\le K,
    u\in\mathcal P(\Gamma)}
  c(u,e)\langle u, [g]\rangle_\alpha=\\
  =\sum_{|u|\le K, u\in\mathcal P(\Gamma)} (\sum_{e\in E\Upsilon}
  \ell(e) c(u,e))\langle u, [g]\rangle_\alpha.
\end{gather*}

Thus the statement of the theorem holds with $d(u)=\sum_{e\in
  E\Upsilon} \ell(e) c(u,e)$.
\end{proof}

\section{Uniform measures and uniform currents}

\begin{defn}[Uniform current]
  Let $A$ be a free basis of $F$ and let $\Gamma$ be a bouquet of $k$
  edges labelled by the generators of $F$ and let $\alpha:F\to
  \pi_1(\Gamma)$ be the corresponding simplicial chart. Thus
  $X=\widetilde \Gamma$ is the Cayley graph of $F$ with respect to
  $A$.
  
  We define a current $n_A$ on $F$ by setting
  $n_A(Cyl(\gamma)):=\frac{1}{2k(2k-1)^{|\gamma|-1}}$ for any path
  $\gamma\in \mathcal P(X)$. It is easy to see that $n_A\in Curr(F)$
  which has weight $1$ with respect to $\alpha$.
  
  We refer to $n_A$ as the \emph{uniform current on $F$ corresponding
    to $A$}.
  
  Similarly, we define a measure $m_A$ on $\partial F$ (identified
  with the set of semi-infinite freely reduced words over $A$) as
  follows. For any freely reduced word $v$ over $A$ we have
  $m_A(Cyl(v)):=\frac{1}{2k(2k-1)^{|v|-1}}$. Then $m_A$ is easily seen
  to be a Borel probability measure on $\partial F$ that is invariant
  with respect to the shift map $T_A:\partial F\to\partial F$ that
  erases the first letter of each geodesic ray. We call $m_A$ the
  \emph{uniform measure on $\partial F$ corresponding to $A$}.
\end{defn}

\begin{conv}
  For a point $\zeta\in \partial F$ and an integer $n\ge 0$ we denote
  by $\zeta_A(n)$ the element of $F$ corresponding to the initial
  segment of $\zeta$ of length $n$, when $\zeta$ is expressed as a
  geodesic ray over $A$.
\end{conv}

The following is an easy corollary of the law of large numbers applied
to the finite state markov process generating freely reduced words
over $A$:
\begin{prop}\label{av}
  Let $A$ be a free basis of $F$. Then for a $m_A$-a.e. point
  $\zeta\in \partial F$ we have
\[
\lim_{n\to\infty} \frac{\eta_{\zeta_A(n)}}{n}=n_A
\]
in $Curr(F)$ and
\[
\lim_{n\to\infty} [\eta_{\zeta_A(n)}]=[n_A]
\]
in $\mathbb PCurr(F)$.
\end{prop}
Informally, the above statement says that in a long random freely
reduced word $\zeta_A(n)$ the frequency of every fixed freely reduced
word $v$ approaches its equilibrium value $\frac{1}{2k(2k-1)^{|v|-1}}$
and that the word $\zeta_A(n)$ is ``almost" cyclically reduced.

\begin{lem}\label{lem:unif}
  Let $A$ be a free basis of $F$ and let $\phi\in Aut(F)$. Then $\phi
  n_A=n_{\phi(A)}$ in $Curr(F)$.
\end{lem}
\begin{proof}
  Let $B=\phi(A)$ and let $v=v(B)$ be a freely reduced word over $B$
  of length $t>0$. Let $u=v(A)$, so that $\phi(u)=v$.
  
  Let $X_A$ and $X_B$ be the Cayley graphs of $F$ with respect to $A$
  and $B$ respectively. Then at the level of subsets of $\partial F$
  we have $\phi(Cyl_A(u))=Cyl_B(v)$. Similarly, at the level of
  subsets of $\partial^2 F$, if $\gamma$ is a path in $X_B$ labelled
  $v=v(B)$ then $\phi^{-1}(Cyl_B(\gamma))=Cyl_A(\beta)$ where $\beta$
  is a path in $X_A$ labelled by $u=v(A)$.
  
  Hence
\[
(\phi
n_A)(Cyl_B(\gamma))=n_A(\phi^{-1}(Cyl_B(\gamma)))=n_A(Cyl_A(\beta))=\frac{1}{2k(2k-1)^{t-1}}.
\]
This implies that $\phi n_A=n_{\phi(A)}$, as claimed.
\end{proof}

The following fact is a basic consequence of the Subadditive Ergodic
Theorem and of the non-amenability of $F$.  It is essentially a
restatement of the results of Kapovich, Kaimanovich and Schupp
in~\cite{KKS} where we refer the reader for a more detailed discussion
about generic stretching factors.

\begin{lede}[Generic stretching factors]
  Let $A$ be a free basis of $F$, let $\ell\in Len(F)$ and let $X$ be
  an $\mathbb R$-tree realizing the length function $\ell$. Then there
  exists a unique number $\lambda_A(\ell)>0$ with the following
  property.
  
  For $m_A$-a.e point $\zeta\in \partial F$ we have:

\begin{enumerate}
\item
  $\displaystyle\lim_{n\to\infty}\frac{\ell(\zeta_A(n))}{n}=\lim_{n\to\infty}\frac{\ell(\zeta_A(n))}{\ell_A(\zeta_A(n))}=\lambda_A(\ell)$;
\item $\displaystyle\lim_{n\to\infty} \frac{d_X(p,
    \zeta_A(n)p)}{n}=\lambda_A(\ell)$ where $p\in X$ is any point.
\end{enumerate}

The number $\lambda_A(\ell)$ is called the \emph{generic stretching
  factor of $\ell$ with respect to $A$}.
\end{lede}
Informally, for a long random freely reduced word $w$ over $A$ we have
$\ell(w)/||w||_A\approx \lambda_A(\ell)$.

\begin{defn}
  Let $F$ be a free group with a free basis $A$ and let $\phi\in
  Out(F)$. We call $\lambda_A(\ell_A\phi)=\lambda_A(\phi^{-1} l_A)$
  the \emph{generic stretching factor of $\phi$ with respect to $A$}
  and denote it by $\lambda_A(\phi)$.
  
  If $\varphi\in Aut(F)$ is an automorphism representing $\phi\in
  Out(F)$, we set $\lambda_A(\varphi):=\lambda_A(\phi)$.
\end{defn}
Thus if $\varphi\in Aut(F)$ represents $\phi$ then for a long random
freely reduced word $w$ over $A$ we have
$\frac{l_A(\varphi(w))}{\ell_A(w)}\approx \lambda_A(\phi)$.

\section{Computation of generic stretching factors}

Let $F=F(A)$ be a free group of finite rank $k\ge 2$ with a free basis
$A=a,b,\dots$. Our goal in this section is to produce an automorphism
$\phi$ of $F$ such that $\lambda_A(\phi)\ne \lambda_A(\phi^{-1})$.

Consider the automorphisms $\tau,\sigma\in Aut(F)$ defined as follows.
We have $\tau(b)=ba$ and $\tau(x)=x$ for each $x\in A-\{b\}$. We have
$\sigma(a)=ab$ and $\sigma(x)=x$ for each $x\in A-\{a\}$. Finally put
$\phi=\sigma\tau^2$. We claim that $\lambda_A(\phi)\ne
\lambda_A(\phi^{-1})$.

We shall need the following series of lemmas for explicit computations
of the generic stretching factors of $\phi$ and $\phi^{-1}$. For a
cyclic word $w$ over $A$ and for a freely reduced word $v$ we will
denote $n(v,w):=\langle v, w\rangle_A+\langle v^{-1}, w\rangle_A$.

\begin{lem} Let $w$ be any nontrivial cyclic word over $A$. Then:
\begin{enumerate}
\item $n(x,\tau(w))=n(x,w)$ for any $x\in A$, $x\ne a$.
  
\item $n(ab^{-1},\tau(w))=n(ba^{-1},
  \tau(w))=n(ba^{-1}b^{-1},w)+n(ba^{-2},w)$.
  
\item $|\tau(w)|=|w|+n(b,w)-2n(ba^{-1},w)=|w|+n(a,\tau(w))-n(a,w)$ and
  therefore $n(a,\tau(w))=n(a,w)+n(b,w)-2n(ba^{-1},w)$.
  
\item $n(ba^{-1}b^{-1}, \tau(w))=n(ba^{-1}b^{-1}, w)$.
\item $n(ba^{-2}, \tau(w))=n(ba^{-3},w)+n(ba^{-2}b^{-1}, w)$.
\end{enumerate}
\end{lem}
\begin{proof}
  Note that $\tau(b^{-1})=a^{-1}b^{-1}$. The lemma follows easily from
  the fact that the only cancellations in $\tau(w)$ after the
  letter-wise application of $\tau$ are of the form $aa^{-1}$ and they
  come from the occurrences of $b a^{-1}$ and $ab^{-1}$ in $w$. In
  particular, no letter different from $a^{\pm 1}$ is cancelled. After
  these cancellations of $aa^{-1}$ are performed, the result is the
  cyclically reduced form of $\tau(w)$.
\end{proof}

The following two lemmas are essentially self-explanatory and we omit
the details:

\begin{lem} Let $w$ be any nontrivial cyclic word over $A$. Then:
\begin{enumerate}
  
\item $n(x,\tau^2(w))=n(x,w)$ for any $x\in A$, $x\ne a$.
  
\item $n(ab^{-1},\tau^2(w))=n(ba^{-1},
  \tau^2(w))=n(ba^{-1}b^{-1},\tau(w))+n(ba^{-2},\tau(w))=
  n(ba^{-1}b^{-1}, w)+n(ba^{-3},w)+n(ba^{-2}b^{-1}, w)$.
  
\item $n(a,\tau^2(w))=n(a,\tau(w))+n(b,\tau(w))-2n(ba^{-1},\tau(w))=
  n(a,w)+n(b,w)-2n(ba^{-1},w)
  +n(b,w)-2[n(ba^{-1}b^{-1},w)+n(ba^{-2},w)]=n(a,w)+2n(b,w)-2n(ba^{-1},w)-2[n(ba^{-1}b^{-1},w)+n(ba^{-2},w)]$.
  
\item
  $|\tau^2(w)|=|\tau(w)|+n(b,\tau(w))-2n(ba^{-1},\tau(w))=|w|+n(b,w)-2n(ba^{-1},w)+n(b,w)-2[n(ba^{-1}b^{-1},w)+n(ba^{-2},w)]=
  |w|+2n(b,w)-2n(ba^{-1},w)-2[n(ba^{-1}b^{-1},w)+n(ba^{-2},w)]$.
\end{enumerate}
\end{lem}

\begin{lem}\label{lem:phi} Let $w$ be any nontrivial cyclic word over $A$. Then:
\begin{enumerate}
  
\item $|\sigma(w)|=|w|+n(a,w)-2n(ab^{-1},w)$
  
\item $|\sigma\tau^2w|=|\tau^2w|+n(a,\tau^2w)-2n(ab^{-1},\tau^2w)=
  |w|+2n(b,w)-2n(ba^{-1},w)-2[n(ba^{-1}b^{-1},w)+n(ba^{-2},w)]+
  n(a,w)+2n(b,w)-2n(ba^{-1},w)-2[n(ba^{-1}b^{-1},w)+n(ba^{-2},w)] -
  2[n(ba^{-1}b^{-1},w)+n(ba^{-3},w)+n(ba^{-2}b^{-1}, w)]=
  |w|+4n(b,w)+n(a,w)-4n(ba^{-1},w)-6n(ba^{-1}b^{-1},w)-4n(ba^{-2},w)-2n(ba^{-3},w)-2n(ba^{-2}b^{-1},w).$
\end{enumerate}
\end{lem}

\begin{prop}\label{phi}
  We have
  $\lambda_A(\phi)=1+\frac{5}{k}-\frac{4}{k(2k-1)}-\frac{10}{k(2k-1)^{2}}-\frac{4}{k(2k-1)^{3}}$.
\end{prop}

\begin{proof}
  If $w$ is a long random cyclic word, then the frequency of a fixed
  reduced word $v$ of length $t$ tends to the uniform frequency
  $\frac{1}{2k(2k-1)^{t-1}}$. Therefore $n(v,w)/|w|$ tends to
  $r(t):=\frac{1}{k(2k-1)^{t-1}}$.
  
  Hence by part (2) of Lemma~\ref{lem:phi} we have:
\begin{gather*}
  \lambda_A(\phi)=1+5r(1)-4r(2)-10r(3)-4r(4)=\\
  1+\frac{5}{k}-\frac{4}{k(2k-1)}-\frac{10}{k(2k-1)^{2}}-\frac{4}{k(2k-1)^{3}}.
\end{gather*}
\end{proof}

We now need to perform similar computations for
$\phi^{-1}=\tau^{-2}\sigma^{-1}$. Note that $\tau^{-1}(b)=ba^{-1}$,
$\sigma^{-1}(a)=ab^{-1}$ and that $\tau^{-1}$ and $\sigma^{-1}$ fix
all other letters.

The following two lemmas are again, easy corollaries of the
definitions:
\begin{lem} Let $w$ be any nontrivial cyclic word over $A$. Then:
\begin{enumerate}
  
\item
  $|\tau^{-1}(w)|=|w|+n(b,w)-2n(ba,w)=|w|+n(a,\tau^{-1}(w))-n(a,w)$
  and hence $n(a,\tau^{-1}(w))=n(a,w)+n(b,w)-2n(ba,w)$.
  
\item $n(ba,\tau^{-1}w)=n(ba^2,w)$
  
\item $|\tau^{-2}(w)|=|\tau^{-1}w|+n(b,\tau^{-1}w)-2n(ba,\tau^{-1}w)=
  |w|+n(b,w)-2n(ba,w)+n(b,w)-2n(ba^2,w)=|w|+2n(b,w)-2n(ba,w)-2n(ba^2,w)$.
\end{enumerate}
\end{lem}

\begin{lem} Let $w$ be any nontrivial cyclic word over $A$. Then:
\begin{enumerate}
  
\item
  $|\sigma^{-1}(w)|=|w|+n(a,w)-2n(ab,w)=|w|+n(b,\sigma^{-1}(w))-n(b,w)$
  and hence $n(b, \sigma^{-1}w)=n(a,w)+n(b,w)-2n(ab,w)$;

\item $n(ba, \sigma^{-1}w)=n(ba,w)-n(aba,w)$;
  
\item $n(ba^2,\sigma^{-1}w)=n(w,baba)-n(ababa,w)$;
  
\item
  $|\tau^{-2}\sigma^{-1}w|=|w|+2n(b,\sigma^{-1}w)-2n(ba,\sigma^{-1}w)-2n(ba^2,\sigma^{-1}w)=
  |w|+2[n(a,w)+n(b,w)-2n(ab,w)]-2n(ba,w)+2n(aba,w)-2n(w,baba)+2n(ababa,w)$.
\end{enumerate}
\end{lem}

\begin{prop}\label{phi-1}
  We have
  \[\lambda_A(\phi^{-1})=1+\frac{4}{k}-\frac{6}{k(2k-1)}+\frac{2}{k(2k-1)^{2}}-\frac{2}{k(2k-1)^{3}}+\frac{2}{k(2k-1)^{4}}.\]
\end{prop}

\begin{proof}
  We will use the same notations as in the proof of
  Proposition~\ref{phi}. Then it follows from the previous lemma that
\begin{gather*}
  \lambda_A(\phi)=1+4r(1)-6r(2)+2r(3)-2r(4)+2r(5)=\\
  1+\frac{4}{k}-\frac{6}{k(2k-1)}+\frac{2}{k(2k-1)^{2}}-\frac{2}{k(2k-1)^{3}}+\frac{2}{k(2k-1)^{4}}.
\end{gather*}
\end{proof}

Proposition~\ref{phi} and Proposition~\ref{phi-1} immediately imply:
\begin{cor}\label{differ}
  Let $F$ be a finitely generated free group of rank $k\ge 2$. Then
  for any free basis $A$ of $F$ there exists an automorphism $\phi$ of
  $F$ such that $\lambda_A(\phi)\ne \lambda_A(\phi^{-1})$.
\end{cor}

\section{Interpreting the intersection form as the distortion of a random geodesic}

Proposition~\ref{av} yields a geometric interpretation of the value of
the intersection form $I(\ell,n_A)$ where $\ell\in FLen(F)$ is
arbitrary and $n_A$ is the uniform current corresponding to a free
basis $A$ of $F$. This is similar to Bonahon's interpretation of the
intersection number between Liouville currents corresponding to two
hyperbolic structures on a compact surface as the generic distortion
of a long random geodesic in first hyperbolic structure with respect
to the second hyperbolic structure.

\begin{prop}\label{interpret}
  Let $A$ be a free basis of $F$ and let $\ell\in FLen(F)$.
  
  Then
\[
I(\ell,n_A)=\lambda_A(\ell).
\]

In particular, for an arbitrary $\phi\in Out(F)$

\[
I(\phi^{-1}\ell_A,n_A)=I(\ell_A, \phi n_A)=I(\ell_A,
n_{\phi(A)})=\lambda_A(\phi).
\]
\end{prop}
\begin{proof}
  Let $\zeta\in \partial F$ be a $m_A$-random point. Then by
  Proposition~\ref{av} we have $\displaystyle\lim_{n\to\infty}
  \frac{\eta_{\zeta_A(n)}}{n}=n_A$. By continuity of $I$ it follows
  that
\[
I(\ell,n_A)=\lim_{n\to\infty} \frac{1}{n} I(\ell,
\eta_{\zeta_A(n)})=\lim_{n\to\infty}
\frac{\ell(\zeta_A(n))}{n}=\lambda_A(\ell).
\]
\end{proof}

We should also note that the generic stretching factor
$\lambda_A(\ell)$ can be interpreted in terms of the Hausdorff
dimension of the measure $m_A$ with respect to the metric on $\partial
F$ corresponding to $F$. We refer the reader to \cite{Kaim,KKS} for
more details.

Recall that for the hyperbolic surface case Bonahon's notion of an
intersection number between two geodesic currents is symmetric. We can
now prove that in the free group case such symmetry is essentially
impossible.

\begin{thm}\label{asym}
  There does not exist a symmetric $Out(F)$-equivariant map
\[
\hat I: Curr(F)\times Curr(F)\to \mathbb R
\] 
such that for some free basis $A$ of $F$, for some $r>0$ and for every
$\nu\in Curr(F)$ we have
\[
\hat I(n_A,\nu)=r I(\ell_A, \nu).
\]
\end{thm}
\begin{proof}
  Suppose, on the contrary, that such a map $\hat I$, a free basis $A$
  of $F$ and anumber $r>0$ exist.

  Let $\phi\in Out(F)$. By $Out(F)$-equivariance and symmetry we have

\begin{gather*}
  I(\ell_A, \phi n_A)=\frac{1}{r}\hat I(n_A, \phi n_A)=\frac{1}{r}\hat
  I(\phi^{-1}n_A,
  n_A)=\text{ by symmetry }\\
  \frac{1}{r}\hat I(n_A, \phi^{-1}n_A)=I(\ell_A, \phi^{-1}n_A).
\end{gather*}
However, by Proposition~\ref{interpret} we have $I(\ell_A, \phi
n_{A})=\lambda_A(\phi)$ and $I(\ell_A, \phi^{-1}
n_{A})=\lambda_A(\phi^{-1})$. Hence for every $\phi\in Aut(F)$ we have
$\lambda_A(\phi)=\lambda_A(\phi^{-1})$, which contradicts
Corollary~\ref{differ}.
\end{proof}

\section{Finite-dimensional approximations}\label{sect:approx}

\begin{notation}
  Let $\alpha: F\to \pi_1(\Gamma,p)$ be a simplicial chart for $F$.
  Let $E\Gamma$ denote the set of oriented edges of $\Gamma$. If $u$
  is an edge-path of positive length in $\Gamma$, we denote by $a(u)$
  the set of all $e\in E\Gamma$ such that $eu$ is a reduced edge-path
  in $\Gamma$. Similarly, denote by $b(u)$ the set of all $e\in
  E\Gamma$ such that $ue$ is a reduced edge-path in $\Gamma$.
\end{notation}

\begin{defn}
  Let $m\ge 1$ be an integer. Let $S(m)=S_\Gamma(m)$ denote set of all
  $v\in \mathcal P(\Gamma)$ with $|v|=m$ and let $D(m)=D_\Gamma(m)$ be
  the number of elements in $S(m)$. We will think of points of
  $\mathbb R^{D(m)}$ as tuples $x=(x_v)_{v\in S(m)}$.
  
  Put

\begin{gather*}
  R_m=R_m(\Gamma):=\{x=(x_v)_{v\in S(m)}\in \mathbb R^{D(m)}:\\
  x_v\ge 0 \text{ for each } v\in S(m), \text{ and}\\
  \sum_{e\in a(u)} x_{eu}=\sum_{e\in b(u)} x_{ue} \text{ for each }
  u\in S_\Gamma(m-1).  \}
\end{gather*}

For a point $x\in R_m(\Gamma)$ denote $\omega(x):=\sum_{v\in S(m)}
x_v$ and call it the \emph{weight} of $x$.

Put $Q_m=Q_m(\Gamma):=\{x\in R_m(\Gamma):\omega(x)=1\}$.
\end{defn}

Thus both $Q_m$ and $R_m$ are finite-dimensional convex polyhedra and,
in addition, $Q_m$ is compact.

\begin{lede}
  Define $\pi_m: \mathbb R^{D(m)}\to \mathbb R^{D(m-1)}$ as follows:

\[
\pi_m:(x_v)_{v\in S(m)}\mapsto (x_u)_{u\in S(m-1)}
\]
where for every $u\in S(m-1)$ $x_u:=\sum_{e\in a(u)} x_{eu}$.

Then:

\begin{enumerate}
\item We have $\pi_m(R_m)=R_{m-1}$ and $\pi_m(Q_m)=Q_{m-1}$.
  
\item We have $\omega(\pi_m(x))=\omega(x)$ for every $x\in R_m$.
\end{enumerate}
\end{lede}

\begin{lem}
  Suppose $\nu\in Curr(F)$. For each $v\in\mathcal P(\Gamma)$ denote
  $x_v:=\langle v, \nu\rangle$.  Then for each $m\ge 1$ the point
  $x=(x_v)_{v\in S(m)}$ belongs to $R_m$ and
  $\omega_\alpha(\nu)=\omega(x)$. In particular, if $\nu$ is
  $\alpha$-normalized then $x\in Q_m$.
\end{lem}
\begin{proof}
  It is clear that all $x_v\ge 0$. Let $u\in S(m-1)$ and let
  $\gamma\in \mathcal P(X)$ be a lift of $u$.
  
  Then
\[
Cyl(\gamma)=\sqcup_{e\in a(\gamma)} Cyl(e\gamma)=\sqcup_{e\in
  b(\gamma)} Cyl(\gamma e)
\]
Since $\nu$ is finitely-additive,
\[
\nu(Cyl(\gamma))=\sum_{e\in a(\gamma)} \nu(Cyl(e\gamma))=\sum_{e\in
  b(\gamma)} \nu(Cyl(\gamma e)),
\]
that is $x_u=\sum_{e\in a(u)} x_{eu}=\sum_{e\in b(u)} x_{ue}$.

Since $u\in S(m-1)$ was arbitrary, this means that $x=(x_v)_{v\in
  S(m)}$ belongs to $R_m$, as claimed.
\end{proof}

\begin{conv}
  We denote by $j_m:Curr(F)\to R_m$ the map that sends each $\nu\in
  Curr(F)$ to $(\langle v, \nu\rangle)_{v\in S(m)}$.  We denote by
  $\overline {j}_m: \mathbb PCurr(F)\to Q_m$ the quotient of the map
  $j_m$ at the level of projectivizations. That is $\overline
  j_m([\nu])=j_m(\nu_\alpha)$.
\end{conv}

The following is essentially proved in \cite{Ka}:

\begin{prop}\label{invlim}
  We have canonical homeomorphisms \newline $Curr(F)\cong
  \underset{\longleftarrow}{\lim} (R_m, \pi_m)$ and $\mathbb
  PCurr(F)\cong \underset{\longleftarrow}{\lim} (Q_m, \pi_m)$.
  
  Moreover, the maps $j_m: Curr(F)\to R_m$ and $\overline j_m: \mathbb
  PCurr(F)\to Q_m$ are ``onto'' for $m\ge 2$.
\end{prop}

Because of the above proposition we think of $R_m$ and $Q_m$ as
finite-dimensional approximations to $Curr(F)$ and $\mathbb PCurr(F)$
respectively.

\begin{notation}
  Let $v\in \mathcal P(\Gamma)$ be a path of length $m\ge 2$, so that
  $x\in S(m)$. We denote by $v-$ the initial segment of $v$ of length
  $m-1$ and we denote by $v+$ the terminal segment of $v$ of length
  $m-1$. Thus $v-, v+ \in S(m-1)$.
\end{notation}

\begin{defn}[Initial graph]
  Let $m\ge 2$ and let $x=(x_v)_{v\in S(m)}$ be a point in
  $R_m(\Gamma)$. The \emph{initial graph $\Delta(x)$ of $x$} is
  defined as follows: $\Delta(x)$ is a directed labelled graph with
  $V\Delta(x):=S(m-1)$ and $E\Delta(x)=S(m)$. For an edge $v\in S(m)$
  of $\Delta(x)$ the initial vertex of $v$ is $v-$ and the terminal
  vertex of $v$ is $v+$. The edge $v$ is labelled by the number $x_v$.
  
  We also define the \emph{improved initial graph $\Delta'(x)$ of $x$}
  as the union of all edges of $\Delta(x)$ with positive labels,
  together with the end-vertices of these edges.

  For a vertex $u$ of $\Delta(x)$ the sum of the labels on the
  incoming edges at $u$ is equal to the sum of the labels on the
  outgoing edges from $u$. We denote this sum by $d_x(u)$.
\end{defn}

Note that for any $x\in S(m)$ we have $\sum_{u\in S(m-1)}
d_x(u)=\omega(x)=\sum_{v\in S(m)} x_v$.

If $g\in F$ is a nontrivial element, then for the point
$x=\alpha_m(\eta_g)$ all coordinates $x_v$ are integers. Namely, $x_v$
is the number of occurrences of $v$ in the cyclic path $c(g)$ in
$\Gamma$ representing $[g]$. It is natural to ask which points of
$R(m)$ with integer coordinates arise in this way. It turns out that
one can provide an explicit answer in terms of initial graphs.

\begin{prop}\label{realize}
  Let $m\ge 2$ and let $x=(x_v)_{v\in S(m)}$ be a nonzero point in
  $R_m(\Gamma)$.  Then there is $g\in F$ with $x=\alpha_m(\eta_g)$ if
  and only if all the coordinates $x_v$ of $v$ are integers and the
  improved initial graph $\Delta'(x)$ of $x$ is topologically
  connected.
\end{prop}

\begin{proof}
  This statement is essentially proved in \cite{Ka} (where the problem
  of which rational points of $Q_m$ are realized as the images of
  normalized rational currents is considered) and the proof is exactly
  the same here.
  
  We will sketch the argument for the ``if" direction. Suppose all
  $x_v$ are integers and $\Delta'(x)$ is connected. Let $N=\omega(x)$.
  Thus $N>0$ is an integers. If we think of $\Delta'(x)$ as a directed
  multi-graph, where each edge with a label $n>0$ is thought of as $n$
  multiple edges, then for each vertex $u$ the out-degree is equal to
  the in-degree at $u$. Since $\Delta'(x)$ is connected, there exists
  an \emph{Euler circuit} in $\Delta'(x)$, that is a circuit $\gamma$
  such that for each directed edge $v$ of $\Delta'(x)$ the circuit
  $\gamma$ passes through the edge $v$ exactly $x_v$ times. Note that
  $|\gamma|=N$. We obtain a cyclic path $c$ in $\Gamma$ from the
  circuit $\gamma$ as follows. We replace each edge $v\in S(m)$ in
  $\gamma$ by the last edge $e$ of $v$, when $v$ is considered as a
  path in $\Gamma$. The result is a reduced cyclic path $c$ in
  $\Gamma$ and, as is easily seen, for every $v\in S_\Gamma(m)$ the
  number of occurrences of $v$ in $c$ is equal to $x_v$. Thus if $g\in
  F$ is an element represented by $c$, then $x=\alpha_m(\eta_g)$, as
  required.
\end{proof}

\begin{prop}\label{extremal}
  Let $x$ be a nonzero point of $Q_m$ for $m\ge 2$. Then $x$ is
  extremal if and only if the improved initial graph $\Delta'(x)$ is
  isomorphic to a directed simplicial circle where the labels of all
  edges are equal.
  
  In particular, $\Delta'(x)$ is connected for extremal points, and
  hence all extremal points of $Q_m$ are $j_m$-images of rational
  currents.
\end{prop}

\begin{proof}

  Recall that $Q_m$ is a convex finite-dimensional compact polyhedron.
  
  Let $x$ be an extremal point of $Q_m$.
  
  Choose a subgraph $\Lambda$ of $\Delta'(x)$ such that $\Lambda$ is a
  directed simplicial circle. Then there exists $y\in R_m$ such that
  $\Delta'(y)$ is $\Lambda$ where each edge of $\Lambda$ is given
  label $1$. We claim that $x$ is a scalar multiple of $y$.  Suppose
  not.

  Then by looking at the definition of $R_m$ we see that for a
  sufficiently small $\epsilon>0$ the points $x+\epsilon y$ and
  $x-\epsilon y$ of $\mathbb R^{D(m)}$ belong to $R_m$. Note that
  $\omega(x)=1=\frac{1}{2}(\omega(x+\epsilon y)+\omega(x-\epsilon
  y))$.  Since some edge of $\Delta'(x)$ is not contained in
  $\Lambda$, the points $x\pm \epsilon y$ are not scalar multiples of
  $x$. Thus $x$ is a convex linear combination of two points of $R_m$
  that are not scalar multiples of $x$:
\[
x=(\omega(x-\epsilon y)/2)\frac{x-\epsilon y}{\omega(x-\epsilon y)}+
(\omega(x+\epsilon y)/2)\frac{x+\epsilon y}{\omega(x+\epsilon y)}.
\]
It follows that $x$ is a convex linear combination of two points of
$Q_m$ different from $x$, contrary to our assumption that $x$ is an
extremal point of $Q_m$. This proves the ``only if'' implication of
the proposition. We leave the ``if'' direction to the reader.
\end{proof}

The above proof is due to a UIUC Geometric Group Theory REU student
Tyler Smith.

We have seen before that there are natural continuous linear maps
$\pi_m: R_m\to R_{m-1}$. It turns out that there are also canonical
continuous (but not linear) maps going in the opposite direction. They
arise from performing certain kinds of random walks in initial graphs.

\begin{defn}\label{walk}
  Let $m\ge 1$ and let $x=(x_v)_{v\in S(m)}$ be a nonzero point in
  $R_m(\Gamma)$. We consider the following random walk on the initial
  graph $\Delta(x)$.
  
  The initial distribution $\theta=\theta_x$ on the vertex set
  $S(m-1)$ of $\Delta(x)$ is given as $\theta(u):=d_x(u)/\omega(x)$.
  
  For two vertices $u, u'\in S(m-1)$ of $\Delta(x)$ the transition
  probability $\rho(u,u')$ is set to be $0$ if there is no directed
  edge from $u$ to $u'$ in $\Delta(x)$ and it is defined as
  $\rho_x(u,u'):=x_v/d_x(u)$ if there is an edge $v\in S(m)$ from $u$
  to $u'$ in $\Delta(x)$.
\end{defn}

Note that for an edge $v\in S(m)$ of $\Delta(x)$ the
$P_\theta$-probability that the trajectory of the random walk begins
with $v-,v+$ is $\theta_x(v-)\rho_x(v-,v+)=d_x(v-)
\frac{x_v}{d_x(v-)}=x_v$.

\begin{lem}\label{y}
  Let $m\ge 1$ and let $x=(x_v)_{v\in S(m)}$ be a nonzero point in
  $R_m(\Gamma)$. Consider the initial distribution $\theta=\theta_x$
  on the vertex set $S(m-1)$ of $\Delta(x)$ and the random walk on
  $\Delta(x)$ as in Definition~\ref{walk}.
  
  For each $z\in S_\Gamma(m+1)$ put
\[
y_z:=x_{z-}\rho_x(z-+, z++)=x_{z-} \frac{x_{z+}}{d_x(z-+)}.
\]
Thus, in view of the above remark, $y_z$ is the $P_\theta$-probability
that the trajectory of the random walk in $\Delta(x)$ begins with
$z--, z-+,z++$, that is, the walk begins by going through the
edge-sequence $z-,z+$ of $\Delta(x)$.

Then the point $y=(y_z)_{z\in S_\Gamma(m+1)}$ belongs to $R_m(\Gamma)$
and $\pi_{m+1}(y)=x$.
\end{lem}

\begin{proof}
  Let $v\in S(m)$. We need to verify that $x_v=\sum_{e\in a(v)}
  y_{ev}=\sum_{h\in b(v)} y_{vh}$.
  
  Then $v$ is an edge from $v-$ to $v+$ in $\Delta(x)$. Let
  $t_1,\dots, t_{l}$ be the labels of the incoming edges for the
  vertex $v-$ of $\Delta(x)$ and let $s_1, \dots, s_{n}$ be the labels
  of the outgoing edges from $v+$ in $\Delta(x)$. Thus $t_1+\dots +
  t_{l}=d_x(v-)$ and $s_1+\dots + t_{n}=d_x(v+)$ since $x\in R_m$.
  
  Then
\[
\sum_{e\in a(v)} y_{ev}=\sum_{i=1}^{l} t_i \frac{x_v}{d_x(v-)}=x_v
\]
and
\[
\sum_{h\in b(v)} y_{vh}=\sum_{i=1}^{n} x_v \frac{s_i}{d_x(v+)}=x_v,
\]
as required.
\end{proof}

\begin{defn}
  For a point $x\in R_m(\Gamma)$ we denote the point $y\in
  R_{m+1}(\Gamma)$ defined in Lemma~\ref{y} by $\iota_m(x)$. Thus, in
  view of Lemma~\ref{y} we have a map $\iota_m:R_m(\Gamma)\to
  R_{m+1}(\Gamma)$ such that $\pi_{m+1}\circ \iota_m=Id_{R_m}$ (that
  is, $\iota_m$ is a section of $\pi_{m+1}$).
\end{defn}

We list some of the basic properties of $\iota_m$ in the following
statement whose proof is left to the reader:

\begin{prop}
  Let $m\ge 1$. Then:
  
  (a) The map $\iota_m:R_m(\Gamma)\to R_{m+1}(\Gamma)$ is continuous.

  (b) Let $x\in R_m(\Gamma)$, $n\ge m$ and let $y=(\iota_{n}\dots
  \circ\iota_m)(x)\in R_{n+1}$. Let $z\in S(n+1)$. Then $y_z$ can be
  computed as follows. For $i=0,\dots, n+1-m$ let $v_i$ be the
  sub-path of $z$ of length $m$ such that $z\equiv z_iv_iz_i'$ where
  $|z_i|=i$. Thus $v_0$ is the initial segment of length $m$ of $z$
  and $v_{n+1-m}$ is the terminal segment of length $m$ of $z$.  Let
  $u_i=v_{i-1}+=v_{i}-$ for $1\le i\le n+1-m$ and let
  $u_0=v_0-,u_{n+2-m}=v_{n+1-m}+$, so that $u_i$ are vertices of
  $\Delta(x)$.
  
  Then

\[
y_z=x_{v_0} \frac{x_{v_1}}{d_x(u_1)}\frac{x_{v_2}}{d_x(u_2)}\dots
\frac{x_{v_{n+1-m}}}{d_x(u_{n+1-m})}.
\]

That is $y_z$ is equal to the $P_{\theta_x}$-probability that the
trajectory of the random walk in $\Delta(x)$ begins with the vertex
sequence
\[
u_0, u_1, \dots,u_{n+2-m}.
\]

(c) Let $m\ge 1$. For each $x\in R_m(\Gamma)$ and $v\in \mathcal
P(\Gamma)$ with $|v|=n$ denote:

\[
x_v=\begin{cases}
  (\iota_{n-1}\dots \circ\iota_m)(x)_v, \quad \text { if } n>m\\
  x_v , \quad \text { if } n=m\\
  (\pi_{n+1}\dots \circ\pi_m)(x)_v, \quad \text { if } n<m.
\end{cases}
\]

Then $\varepsilon_m(x):=(x_v)_{v\in P(\Gamma)}\in
\underset{\longleftarrow}{\lim} (R_n, \pi_n)$ and $\varepsilon_m:
R_m(\Gamma)\to \underset{\longleftarrow}{\lim} (R_n, \pi_n)$ is a
topological embedding that provides a section for the map
$\alpha_m:\underset{\longleftarrow}{\lim} (R_n, \pi_n)\to R_m$, so
that $\varepsilon_m\circ\iota=Id_{R_m}$.
\end{prop}

\section{The distortion functional}\label{sect:distortion}

\begin{defn}\label{dist}[The distortion functional]
  Let $\alpha:F\to \pi_1(\Gamma, p)$ be a simplicial chart and let
  $X=\widetilde \Gamma$. We equip $X$ and $\Gamma$ with simplicial
  metrics, where every edge has length $1$. Let $\ell=\ell_\alpha\in
  FLen(F)$ be the length function corresponding to the action of $F$
  on $X$ via $\alpha$. Let $\ell'\in FLen(F)$ be arbitrary.

  Put
  \[\delta(\nu)=\delta_{\ell',\ell}(\nu):=I(\ell',\nu)/I(\ell,\nu)\]
  for any nonzero current $\nu\in Curr(F)$.
  
  Note that $\delta(\nu)=\delta(r \nu)$ for any nonzero scalar $r$.
  Thus we can define
  $\delta_{\ell',\ell}([\nu]):=\delta_{\ell',\ell}(\nu)$ for any
  $\nu\in [\nu]$ so that $\delta_{\ell',\ell}$ is now defined on
  $\mathbb PCurr(F)$.
\end{defn}

Recall that $I(\ell',\eta_g)=\ell'(g)$ for every nontrivial $g\in F$.
Thus for a nontrivial $g\in F$ we have
\[
\delta_{\ell',\ell}(\eta_g)=I(\ell',\eta_g)/I(\ell,\eta_g)=\ell'(g)/\ell(g).
\]
Recall also that the $\alpha$-normalized representative $\nu_\alpha$
of $[\nu]$ has the property that $\omega(\nu_\alpha)=I(\ell,
\nu_\alpha)=1$ and hence $\delta(\nu_\alpha)=I(\ell', \nu_\alpha)$.
Since $\mathbb PCurr(F)$ can be thought of (via normalization) as a
convex subset of $Curr(F)$, this implies:

\begin{lem}
  In the notations of Definition~\ref{dist} we have that $\delta$ is a
  continuous linear functional on $\mathbb PCurr(F)$.
\end{lem}

The following is an important corollary of the local formulas
established earlier:

\begin{prop}
  In the notations of Definition~\ref{dist} the maximum and the
  minimum values of $\delta:\mathbb PCurr(F)\to \mathbb R$ are
  achieved and they are realized by rational currents.
\end{prop}

\begin{proof}
  The local formulas obtained in Theorem~\ref{thm:local} imply that
  there exists $m\ge 1$ and a linear functional $\bar \delta:
  Q_m(\Gamma)\to \mathbb R$ such that $\bar\delta(\overline
  j_m(\nu))=\delta(\nu_\alpha)$ for every nonzero $\nu$.  Since
  $\overline j_m$ is ``onto'' $Q_m$, it follows that $\sup \delta=\sup
  \bar \delta$ and $\inf \delta=\inf\bar\delta$.  Since $Q_m$ is
  finite-dimensional convex polyhedron, the linear functional
  $\bar\delta$ does achieve its extremal values on $Q_m$ and they are
  attained at extremal points of $Q_m$. By Proposition~\ref{extremal}
  the extremal points of $Q_m$ have connected improved initial graphs
  and hence correspond to the $\overline j_m$-images of rational
  currents. Thus the result follows.
\end{proof}

\section{Currents and monomorphisms}

Recall that the action of $Aut(F)$ on $Curr(F)$ was defined as
follows. For any $\phi\in Aut(F)$, for any $\nu\in Curr(F)$ and for
any $S\subseteq \partial^2 F$ we have $(\phi\nu)(S):=\nu(\phi^{-1}S)$.
We have verified earlier that this indeed defines an action on
$Curr(F)$ and that, moreover, by Proposition~\ref{good}, for every
nontrivial $g\in F$ we have $\phi(\eta_g)=\eta_{\phi(g)}$.

Suppose $F'$ is another finitely generated free group and that
$\phi:F'\to F$ is an injective homomorphism. Then $\phi$ is a
quasi-isometric embedding and hence it defines an equivariant
topological embedding $\hat \phi:\partial^2 F'\to \partial^2 F$. It is
therefore natural to want to emulate the automorphism case here and to
define a map $\phi_\ast: Curr(F')\to Curr(F)$ in a similar way. Thus
in \cite{Ma} R.Martin claims to define a map $\phi_\ast: Curr(F')\to
Curr(F)$ by the same formula as in the automorphism case: for any
$\nu\in Curr(F')$ and for any Borel $S\subseteq \partial^2 F$

\[(\phi\nu)(S):=\nu(\hat\phi^{-1}S)\tag{$\dag$}
\]
It is then claimed in \cite{Ma} (Lemma 13, Section 5.5) that this
defines an injective proper map $\phi_\ast: Curr(F')\to Curr(F)$.

Unfortunately, this approach is quite incorrect, as are the statements
of Lemma 13 and Lemma 14 in Section 5.5 of \cite{Ma}. The problem is
that that the measure $\phi\nu$ on $\partial^2 F$ defined by $(\dag)$,
will not, generally, be $F$-invariant unless the map $\phi$ is ``onto"
(that is, unless $\phi$ is an isomorphism). For example, suppose
$F=F'\ast F''$ and $\phi$ is the inclusion of $F'$ to $F$. Let
$\nu=\eta_g$ for $g\in F'$. Let $h\in F''$ be a nontrivial element and
let $S=\{(g^{-\infty}, g^{\infty})\}$. Then $hS\cap\hat\phi(\partial^2
F')=\varnothing$ and $\hat\phi^{-1}(hS)=\varnothing$. Thus
$(\phi_\ast\eta_g)(S)=1$ while $(\phi_\ast\eta_g)(hS)=0$, so that
$\phi_\ast\eta_g$ is not $F$-invariant and so does not belong to
$Curr(F)$.

Nevertheless, we show that the following approach works for
monomorphisms and that one should generalize Proposition~\ref{good}
rather than the definition of the action of $Aut(F)$ on $Curr(F)$:

\begin{propdfn}[Maps determined by monomorphisms]
  Let $F'$, $F$ be finitely generated nonabelain free groups and let
  $\phi:F'\to F$ be an injective homomorphism.
  
  Then there exists a unique continuous linear map
  $\phi_\ast:Curr(F')\to Curr(F)$ such that for every nontrivial $g\in
  F'$ we have
\[
\phi_\ast(\eta_g)=\eta_{\phi(g)}.
\]
\end{propdfn}
\begin{proof}
  We will sketch the argument and leave the details to the reader.
  
  Note that if such $\phi_\ast$ exists then it is unique since
  rational currents are dense in $Curr(F')$.
  
  Choose a simplicial chart $\alpha: F\to \pi_1(\Gamma,p)$ for $F$ and
  let $X=\widetilde\Gamma$.
  
  Consider the action of $F'$ on $X$ via $\phi$.
  
  There exists a unique minimal $F'$-tree $Y$ such that the quotient
  graph $\Delta:=Y/F'$ provides a simplicial chart $\beta$ for $F'$ in
  the obvious way. Note that each edge of $\Delta$ is ``labelled'' by
  some path in $\Gamma$ and that this labelling is ``folded'' in the
  obvious way.
  
  Hence the same argument as in the proof of Proposition~\ref{local}
  implies that for each $v\in \mathcal P\Gamma$ there are some
  integers $c(u)=c(u,v)\ge 0$, where $|u|\le K$ and $u\in \mathcal
  P(\Delta)$ such that for every cyclic path $w'$ in $\Delta$ defining
  a cyclic path $w$ in $\Gamma$ we have

\[
\langle v,w\rangle_\alpha=\sum_{|u|\le K, u\in \mathcal P\Gamma}
c(u,v) \langle u,w'\rangle_\beta.
\]
This means that for every nontrivial $g\in F'$ we have
\[
\langle v, \phi(g)\rangle_\alpha=\sum_{|u|\le K, u\in \mathcal
  P\Gamma} c(u,v) \langle u,g\rangle_\beta.
\]
We now define a map $\phi_\ast: Curr(F')\to Curr(F)$ by the following
formula: if $\nu\in Curr(F)$ then for every $v\in \mathcal P\Gamma$
\[
\langle v, \phi_\ast\nu\rangle_\alpha:=\sum_{|u|\le K, u\in \mathcal
  P\Gamma} c(u,v) \langle u,\nu\rangle_\beta.\tag{!}
\]
If these formulas indeed define a geodesic current on $F$ then the map
$\phi_\ast$ is continuous and has the property that for every
nontrivial $g\in F'$ we have $\phi_\ast(\eta_g)=\eta_{\phi(g)}$.

Thus it remains to check that $\phi_\ast\nu$ is indeed a geodesic
current on $F$. Put $x_v=\langle v, \phi_\ast\nu\rangle_\alpha$ for
every $v\in \mathcal P\Gamma$.

By Lemma~\ref{invlim} it suffices to show that the infinite tuple
$(x_v)_v$ defines an element of $\underset{\longleftarrow}{\lim}
(R_m(\Gamma), \pi_m)$. That is we need to show that

(1) for every $m\ge 1$ the tuple $x_m:=(x_v)_{|v|=m}$ is an element of
$R_m(\Gamma)$

and

(2) for every $m\ge 2$ we have $\pi_m(x_m)=x_{m-1}$.

In view of (!) both (1) and (2) reduce to verifying that some
explicitly defined continuous linear functions on $Curr(F)$ are equal
to each other. Since we do have that
$\phi_\ast(\eta_g)=\eta_{\phi(g)}$ for every nontrivial $g\in F'$,
this equalities hold on a dense subset of $Curr(F')$ and hence on the
entire $Curr(F')$ as well.
\end{proof}

However, in general one cannot expect the map $\phi_\ast$ to be
injective since the map from the set of conjugacy classes of $F'$ to
the set of conjugacy classes of $F$ induced by $\phi$ need not be
injective. For example, let $F'=F=F(a,b)$ and let $\phi(a)=a,
\phi(b)=bab^{-1}$. Then $\phi$ is injective but $\phi(a)$ is conjugate
to $\phi(b)$ and hence $\phi_\ast(\eta_a)=\phi_\ast(\eta_b)$.

\section{Translation Equivalence}

The following notion was introduced and studied in detail by Kapovich,
Levitt, Schupp and Shpilrain in \cite{KLSS}.

\begin{defn}[Translation Equivalence of Elements]
  Elements $g,h\in F$ are said to be \emph{translation equivalent in
    F}, denoted $g\equiv_t h$, if for every $\ell\in FLen(F)$ we have
  \[
  \ell(g)=\ell(h).
  \]
\end{defn}
In view of the properties of the intersection form $I$ it makes sense
to generalize it as follows:

 \begin{defn}[Translation Equivalence of Currents]
   Currents $\nu_1,\nu_2\in Curr(F)$ are said to be \emph{translation
     equivalent in $Curr(F)$}, denoted $\nu_1\equiv_t \nu_2$, if for
   every $\ell\in FLen(F)$ we have
 \[
 I(\ell, \nu_1)=I(\ell, \nu_2).
 \]
 \end{defn}
 
 Thus for nontrivial elements $g,h\in F$ we have $g\equiv_t h$ in $F$
 iff $\eta_g\equiv_t \eta_h$ in $Curr(F)$.  The notion of translation
 equivalence can be thought of as measuring the ``degeneracy" of the
 intersection form $I$.

 The following is proved in \cite{KLSS} (recall that by convention $F$
 denotes a free group of finite rank $k\ge 2$.)

\begin{prop}\label{palind}
  Let $F(a,b)$ be free of rank two nd let $\vartheta:F(a,b)\to F(a,b)$
  be the automorphism defined as $\vartheta(a)=a^{-1}$ and
  $\vartheta(b)=b^{-1}$. Then for every $w\in F(a,b)$ and for any
  homomorphism $\phi:F(a,b)\to F$ we have $\phi(w)\equiv_t
  \phi(\vartheta(w))$ in $F$.
\end{prop}
In other words, for any $w(a,b)$ and for any $g,h\in F$ we have
$w(g,h)\equiv_t w(g^{-1}, h^{-1})$ in $F$.  Note that $\vartheta$
represents the only nontrivial central element in $Out(F(a,b)$.

\begin{thm}
  Let $\phi:F(a,b)\to F$ be any injective homomorphism.  Let
  $\vartheta \in Aut(F(a,b))$ be defined as $\vartheta(a)=a^{-1}$ and
  $\vartheta(b)=b^{-1}$.  Then for every $\nu\in Curr(F(a,b)$ we have
\[
\phi_\ast(\nu)\equiv_t \phi_\ast(\vartheta \nu) \text{ in } Curr(F).
\]
\end{thm}

\begin{proof}
  Let $\nu_n\in Curr(F(a,b))$ be a sequence of rational currents such
  that $\displaystyle\lim_{n\to\infty} \nu_n=\nu$. Then
  $\displaystyle\lim_{n\to\infty} \phi_\ast \nu_n=\phi_\ast(\nu)$ and
  $\displaystyle\lim_{n\to\infty} \phi_\ast
  \vartheta\nu_n=\phi_\ast(\vartheta\nu)$.  Let $\mathbb \ell\in
  FLen(F)$ be arbitrary. By Proposition~\ref{palind} $\phi_\ast
  \nu_n\equiv_t \phi_\ast \vartheta\nu_n$ in $Curr(F)$.  Therefore
\begin{gather*}
  I(\ell,\phi_\ast(\nu))=\lim_{n\to\infty} I(\ell,
  \phi_\ast(\nu_n))=\lim_{n\to\infty} I(\ell,
  \phi_\ast(\vartheta\nu_n))=I(\ell, \phi_\ast(\vartheta \nu)),
\end{gather*}
as required.
\end{proof}

\end{document}